%% file: 872.tex
\documentclass[letterpaper]{amsart}
\usepackage{amsfonts,amssymb,amscd,amsmath,amsthm,txfonts}
\usepackage{ifpdf}
\usepackage{accents}
\usepackage[all]{xy}
\ifpdf
   \usepackage[pdftex]{hyperref}
   \usepackage[pdftex]{graphicx}
   \usepackage{epstopdf}
\else
   \usepackage[dvips]{hyperref}
   \usepackage[dvips]{graphicx} 
\fi
\usepackage{color}

\newtheorem{MainThm}{Theorem}
\newtheorem{Thm}{Theorem}[section]
\newtheorem{Lem}[Thm]{Lemma}

\newtheorem{Cor}[Thm]{Corollary}
\theoremstyle{definition}
\newtheorem{Def}[Thm]{Definition}
\newtheorem{Asm}[Thm]{Assumption}
\theoremstyle{remark}

\newtheorem{Note}[Thm]{Note}
\newtheorem*{uRem}{Remark}
\newtheorem*{uRems}{Remarks}

\newtheorem{Facts}[Thm]{Facts}

\numberwithin{equation}{section}

\DeclareMathOperator{\EXP}{EXP}
\DeclareMathOperator{\flarge}{rapidgrowth}
\DeclareMathOperator{\dom}{dom}
\DeclareMathOperator{\myinv}{char}
\newcommand{\myinvset}{\text{CHARS}}
\newcommand{\fnto}{\ensuremath{\rightarrow}}  
\newcommand{\card}[1]{\vert #1 \vert }     
\newcommand{\al}[1]{\ensuremath{{\aleph_{#1}}} }          
\newcommand{\om}[1]{\ensuremath{{\omega_{#1}}} }          
\newcommand{\ho}{\ensuremath{^{\omega}}}                  
\newcommand{\DEFEQ}{\ensuremath{\coloneqq}}
\newcommand{\EQDEF}{\ensuremath{\eqqcolon}}

\newcommand{\SlFam}{\ensuremath{{\mathcal S}}}

\newcommand{\esm}{\ensuremath{\prec}}  

\newcommand{\std}[1]{\ensuremath{\check{#1}}}              
\newcommand{\forc}{\ensuremath{\Vdash}}
\newcommand{\incomp}{\ensuremath{\perp}}
\newcommand{\comp}{\ensuremath{\parallel}}
\newcommand{\n}[1]{\underaccent{\tilde}{#1}}

\DeclareMathOperator{\image}{image}

\DeclareMathOperator{\nor}{nor}
\DeclareMathOperator{\val}{val}
\newcommand{\prodval}{\val^{\Pi}}
\DeclareMathOperator{\chalf}{half}
\DeclareMathOperator{\ns}{preprenor}
\DeclareMathOperator{\prenor}{prenor}
\newcommand{\cc}{\mathfrak{c}}
\newcommand{\cd}{\mathfrak{d}}
\newcommand{\cS}{\mathbf{\Sigma}}
\newcommand{\cSc}{\cS(\cc)}
\newcommand{\cK}{\mathbf{K}}
\newcommand{\bfH}{\mathbf{H}}
\newcommand{\bfF}{\mathbf{F}}
\newcommand{\pow}{\mathcal{P}}
\DeclareMathOperator{\trunk}{trunk}
\DeclareMathOperator{\trunklg}{trnklh}
\DeclareMathOperator{\supp}{supp}

\newcommand{\Qsu}{\mathbb{Q}^*_\infty}
\newcommand{\Qswu}{\mathbb{Q}^*_{w\infty}}
\newcommand{\Qsf}{\mathbb{Q}^*_{f}}
\newcommand{\myc}{c^{\exists}}
\newcommand{\mycfa}{c^{\forall}}

\begin{document}
\date{\today}

\title{Decisive creatures and large continuum}
\author[Jakob Kellner]{Jakob Kellner$^*$}
\address{Kurt G\"odel Research Center for Mathematical Logic\\
 Universit\"at Wien\\
 W\"ahringer Stra\ss e 25\\
 1090 Wien, Austria}
\email{kellner@fsmat.at}
\urladdr{http://www.logic.univie.ac.at/$\sim$kellner}
\thanks{$^*$ supported by a European Union Marie Curie EIF Fellowship,
contract MEIF-CT-2006-024483.}
\author[Saharon Shelah]{Saharon Shelah$^\dag$}

\address{Einstein Institute of Mathematics\\
Edmond J. Safra Campus, Givat Ram\\
The Hebrew University of Jerusalem\\
Jerusalem, 91904, Israel\\
and
Department of Mathematics\\
Rutgers University\\
New Brunswick, NJ 08854, USA}
\email{shelah@math.huji.ac.il}
\urladdr{http://shelah.logic.at/}
\thanks{
$^\dag$
supported by the United States-Israel
  Binational Science Foundation (Grant no. 2002323), and by
the US National Science Foundation grant NSF-DMS 0600940,
publication 872.}


\begin{abstract}
  For $f,g\in\omega\ho$ let $\mycfa_{f,g}$ be the minimal
  number of uniform $g$-splitting trees needed to
  cover the uniform $f$-splitting tree, i.e.\ for every
  branch $\nu$ of the $f$-tree, one of the 
  $g$-trees contains $\nu$.
  $\myc_{f,g}$ is the dual notion: For every branch $\nu$,
  one of the $g$-trees guesses $\nu(m)$ infinitely often.

  It is consistent that
  $\myc_{f_\epsilon,g_\epsilon}=\mycfa_{f_\epsilon,g_\epsilon}=\kappa_\epsilon$
  for $\al1$ many pairwise different cardinals $\kappa_\epsilon$
  and suitable pairs $(f_\epsilon,g_\epsilon)$.

  For the proof we use creatures with
  sufficient bigness and halving. We show that 
  the lim-inf creature forcing satisfies 
  fusion and pure decision. We 
  introduce 
  decisiveness and use it to construct
  a variant of the countable support iteration
  of such forcings, which still satisfies 
  fusion and pure decision.
\end{abstract}
\maketitle

\section{Introduction}
  In the paper {\em Many simple cardinal invariants}~\cite{GoSh:448},
  Goldstern and the second author construct a
  partial order $P$ that forces pairwise different values to
  $\al1$ many instances of
  the cardinal characteristic $\mycfa_{f,g}$,
  defined as follows:

  Let $f,g\in \omega\ho$ (usually we have $f(n)>g(n)$ for all $n$).
  An $(f,g)$-slalom is a sequence
  $S=(S(n))_{n\in\omega}$
  such that $S(n)\subseteq f(n)$ and $\card{S(n)}\leq g(n)$.
  A family $\SlFam$ of
  $(f,g)$-slaloms is a $(\forall,f,g)$-cover, if
  for all $r\in \prod_{n\in\omega} f(n)$ 
  there is an $S\in \SlFam$ such that
  $r(n)\in S(n)$ for all $n\in\omega$.
  $\mycfa_{f,g}$ is the minimal size of a $(\forall,f,g)$-cover.

  We investigate the dual notion:
  A family $\SlFam$ of
  $(f,g)$-slaloms is an $(\exists,f,g)$-cover, if
  for all $r\in \prod_{n\in\omega} f(n)$ there is an $S\in \SlFam$ such that
  $r(n)\in S(n)$ for infinitely many $n\in\omega$.
  $\myc_{f,g}$ is the minimal size of
  an $(\exists,f,g)$-cover.


  In~\cite{GoSh:448}, the following is shown:
  \begin{quote}
    Assume that CH holds, that $(f_\epsilon,g_\epsilon)_{\epsilon\in\om1}$
    are sufficiently different, and that
    \mbox{$\kappa_\epsilon^{\al0}=\kappa_\epsilon$}
    for all $\epsilon\in\om1$. Then there is a cardinal preserving partial
    order $P$ which forces that $\mycfa_{f_\epsilon,g_\epsilon}=\kappa_\epsilon$ for all $\epsilon\in\om1$.
  \end{quote}

  Similar results regarding $\myc$ as well as a perfect set of invariants were
  promised to appear in a paper called 448a, which never materialized.
  A result for continuum many different invariants of the form
  $\mycfa_{f_\epsilon,g_\epsilon}$ can be found in~\cite{morecardinals}.

  In this paper, we prove a version for countably many
  invariants $\myc$:
  \begin{MainThm}\label{thm:ctbl}
    Assume that CH holds, that $(f_\epsilon,g_\epsilon)_{\epsilon\in\omega}$
    are sufficiently different, and that
    \mbox{$\kappa_\epsilon^{\al0}=\kappa_\epsilon$}
    for all $\epsilon\in\omega$. Then there is a cardinal preserving, 
    $\omega^\omega$-bounding partial order $P$ which forces that
    $\myc_{f_\epsilon,g_\epsilon}=\mycfa_{f_\epsilon,g_\epsilon}=\kappa_\epsilon$ for all $\epsilon\in\omega$.
  \end{MainThm}
  (See Section~\ref{sec:countable} for a definition of sufficiently different.)

  We can also get $\om1$ many different invariants, but we do not know
  in the ground model which invariants will be picked:
  \begin{MainThm}\label{thm:uncountable}
    Assume that CH holds, and that
    $\kappa_\epsilon^{\al0}=\kappa_\epsilon$
    for all $\epsilon\in\om1$. Then
    there are pairs $(f_\nu,g_\nu)_{\nu\in \om1}$
    and there is a cardinal preserving,
    $\omega^\omega$-bounding partial
    order $R$ which forces:
    For each $\epsilon\in\om1$ there is a $\nu(\epsilon)\in \om1$
    such that
    $\myc_{f_{\nu(\epsilon)},g_{\nu(\epsilon)}}=
    \mycfa_{f_{\nu(\epsilon)},g_{\nu(\epsilon)}}=\kappa_\epsilon$.
  \end{MainThm}

  In any case, if the $\kappa_\epsilon$ are pairwise different,
  then in the forcing extension there 
  are infinitely many different cardinals below the continuum, i.e.\
  $2^\al0>\al{\omega}$. Therefore we cannot use countable
  support iterations.
  We cannot use finite support iterations either
  (otherwise we add many  Cohen reals, which
  makes $\mycfa$ too big).
  Instead, we use a
  variant of the countable support product of lim-inf
  creature forcings. We do not assume that the reader knows
  anything about creature forcing.
  However, we do assume that the reader knows the definition of proper forcing
  (see e.g.~\cite{MR1234283} or, for the brave,~\cite{MR1623206}),
  and the fact that such forcings preserve $\om1$.
  Alternatively, it is sufficient to know
  Baumgartner's Axiom~A (cf.~\cite{MR823775}): it is easy to see that
  the forcings in this paper all satisfy Axiom~A, and Axiom~A
  forcings (are proper and therefore) preserve $\om1$.

  We write $q\leq p$ to say that $q$ is stronger than $p$.
  We try to stick to Goldstern's alphabetic convention, i.e.\ 
  whenever two conditions are compatible, the symbol used for the
  stronger condition comes lexicographically later.

  The theorems in this paper are due to the second author.
  The first author's contribution was to fill in some details,
  to ask the second author to fill in other details, and 
  to write the paper.

  We thank a referee for very carefully reading the paper and
  pointing out a mistake and numerous unclarities.

  \subsection*{Annotated contents}
  In the first part, we investigate lim-inf creature forcings:
  \begin{list}{}{\setlength{\leftmargin}{0.5cm}\addtolength{\leftmargin}{\labelwidth}}
    \item[Section~\ref{sec:simple}, p. \pageref{sec:simple}.]
      We define the (one-dimensional) 
      lim-inf creature forcing $\Qsu$.
    \item[Section~\ref{sec:puredecondim}, p. \pageref{sec:puredecondim}.]
      We use bigness and halving to show that $\Qsu$ satisfies 
      pure decision (and fusion).
      This implies that  $\Qsu$ is proper and $\omega^\omega$-bounding.
      We also show rapid reading of
      certain names. The proofs in this section will 
      be generalized in Section~\ref{sec:product}.
    \item[Section~\ref{sec:decisive}, p. \pageref{sec:decisive}.]
      We introduce decisiveness and use it to extend bigness to
      functions defined on finite products of creatures. This allows
      us to show pure decision for finite products of lim-inf creature forcings.
    \item[Section~\ref{sec:product}, p. \pageref{sec:product}.]
      We define the forcing $P$, a variant of the countable support product of
      lim-inf creature forcings, in such a way that the proof of
      Section~\ref{sec:puredecondim} still works with only few changes.
      We also get $\al2$-cc (assuming CH).
    \item[Section~\ref{sec:example}, p. \pageref{sec:example}.]
      We show how to construct decisive
      creatures with sufficient bigness and halving.
  \end{list}
  In the second part, we use the methods of Section~\ref{sec:product}
  to prove Theorems~\ref{thm:ctbl} and~\ref{thm:uncountable}:
  \begin{list}{}{\setlength{\leftmargin}{0.5cm}\addtolength{\leftmargin}{\labelwidth}}
    \item[Section~\ref{sec:countable}, p. \pageref{sec:countable}.]
      We formulate the requirements for Theorem~\ref{thm:ctbl}
      and define $P$, a variant the forcing in Section~\ref{sec:product}.
    \item[Section~\ref{sec:forall}, p. \pageref{sec:forall}.]
      We show that $P_\epsilon$, a complete subforcing
      of $P$, adds a $\mycfa_{f_\epsilon,g_\epsilon}$-cover  in $V[G_P]$.
      This proves $\mycfa_{f_\epsilon,g_\epsilon}\leq \kappa_\epsilon$.
    \item[Section~\ref{sec:exists}, p. \pageref{sec:exists}.]
      We show that in $V[G_P]$ there can be no
      $\myc_{f_\epsilon,g_\epsilon}$-cover smaller than $\kappa_\epsilon$:
      Otherwise we can find a condition $q$ that rapidly reads (without
      using index $\beta$)
      a slalom $\n S$ and forces that the generic real
      $\n\eta_\beta$ at $\beta$ 
      meets $\n S$ infinitely often. We strengthen $q$ such that
      the possible values for the generic always\footnote{This
       is the reason we have to use lim-inf creature forcing
       instead of lim-sup: When we deal with $\mycfa$,
       we have to ``run away'' from $\n S$ infinitely often, 
       and it is enough to assume that we have sufficient space to
       do so infinitely often. But here we need sufficient space at {\em every}
       height.}
      avoid the slalom $\n S$, a contradiction.
    \item[Section~\ref{sec:uncountable}, p. \pageref{sec:uncountable}.]
      We construct $\om1$ many suitable pairs $(f_\epsilon,g_\epsilon)$
      the partial order $R$, a modification of $P$, to show
      Theorem~\ref{thm:uncountable}.
  \end{list}

\section{lim-inf creature forcings}\label{sec:simple}
Creature forcing in general is described in the monograph
{\em Norms on possibilities I: forcing with trees and creatures\/}~\cite{RoSh:470} by
Ros\l anowski and the second author. 
The forcing of the proof in~\cite{GoSh:448} can be interpreted as
creature forcing as well, more specifically as a lim-sup tree creating
creature forcing.
We will use lim-inf creatures instead.
These forcings are generally
more complicated than the lim-sup case, and~\cite{RoSh:470}
shows that they can collapse $\om1$.
In this paper, we will require increasingly strong bigness and halving, which
guarantees pure decision and therefore properness.

We now describe the setting we use.
Creature forcings are defined by a parameter, the creating pair $(\cK,\cS)$.
We use the following objects:
\begin{itemize}
  \item A function $\bfH:\omega\fnto\omega\setminus \{0\}$.
  \item A strictly increasing function $\bfF :\omega\fnto \omega$ such that
    $\bfF (0)=0$.
  \item For every $n\in\omega$ a finite set $\cK(n)$.
  \item For each $\cc\in\cK(n)$, a real number
    $\nor(\cc)\geq 0$, and a nonempty subset $\val(\cc)$
    of $\prod_{\bfF (n)\leq i<\bfF (n+1)}\bfH(i)$.
  \item We additionally require that $\card{\val(\cc)}=1$ implies $\nor(\cc)=0$.
\end{itemize}

A $\cc\in\cK(n)$ is called $n$-creature. The intended meaning of the
$n$-creature $\cc$ is the following: the set of possible values for the generic
object $\n\eta\in \prod_{i\in\omega}\bfH(i)$ restricted to the interval
\mbox{$[\bfF (n),\bfF (n+1)-1]$} is the set $\val(\cc)$.  $\nor(\cc)$ can be
thought of measuring the amount of ``freedom'' the creature $\cc$ leaves on its
interval.  If $\cc$ determines its part of the generic real (i.e.\ if
$\val(\cc)$ is a singleton) then $\nor(\cc)=0$ (i.e.\ $\cc$ leaves no freedom).
However, this intuition about $\nor(\cc)$ has to be used with caution:
In particular, $\val(\cd)\subseteq \val(\cc)$ does generally not imply
$\nor(\cd)\leq \nor(\cc)$.

We set $\cK\DEFEQ \bigcup_{n\in\omega} \cK(n)$.

In our application we will use $\bfF (n)=n$, i.e.\ an $n$-creature lives on the
singleton $\{n\}$.

We also have a
function $\cS: \cK\fnto \pow(\cK)$ satisfying:
\begin{itemize}
  \item If $\cc\in\cK(n)$ and $\cd\in\cSc$ then $\cd\in\cK(n)$.
  \item $\cS$ is reflexive, i.e.\ $\cc\in\cSc$.
  \item $\cS$ is transitive, i.e.\ $\cd\in\cSc$ and $\cd'\in\cS(\cd)$ implies
    $\cd'\in\cSc$.
  \item If $\cd\in\cSc$ then $\val(\cd)\subseteq \val(\cc)$ and
    $\nor(\cd)\leq \nor(\cc)$.
\end{itemize}

The intended meaning is that $\cS(\cc)$ is the set of creatures that are
stronger than $\cc$.

To simplify notation later on, 
we extend the definitions of $\nor$, $\val$ and $\cS$ to
sequences $s,t\in\prod_{\bfF (n)\leq i<\bfF (n+1)}\bfH(n)$:
We set
\[
  \nor(t)\DEFEQ  0,\quad
  \val(t)\DEFEQ \{t\},\quad
  t\in \cS(\cc)\text{ iff }t\in\val(\cc),\quad
  s\in \cS(t)\text{ iff }s=t.
\]

We now define the lim-inf forcing $\Qsu(\cK,\cS)$:
\begin{Def}
  A condition $p\in \Qsu(\cK,\cS)$ consists of a trunk
  $t\in \prod_{i<\bfF (n)} \bfH(i)$ for some $n$ and a
  sequence $(\cc_i)_{i\geq n}$
  such that 
  $\cc_i\in\cK(i)$ and $\nor(\cc_i)>0$ for all $i\geq n$,
  and $\lim(\nor(\cc_i))=\infty$.
  We set $\trunk(p)\DEFEQ t$,
  and the trunk-length $\trunklg(p)\DEFEQ n$, and we set
  \[
    p(i)\DEFEQ
    \begin{cases}
      \cc_i&\text{if }i\geq n,\\
      t\restriction[\bfF (i),\bfF (i+1)-1]&\text{otherwise.}
    \end{cases}
  \]
  So we can identify $p$ with the sequence $(p(i))_{i\in\omega}$.
  The order on $\Qsu$ is defined by
  $q\leq p$ if $\trunklg(q)\geq \trunklg(p)$ and
  $q(i)\in \cS(p(i))$ for all $i$.
\end{Def}
\begin{figure}[tb]
  \begin{center}
    \scalebox{0.4}{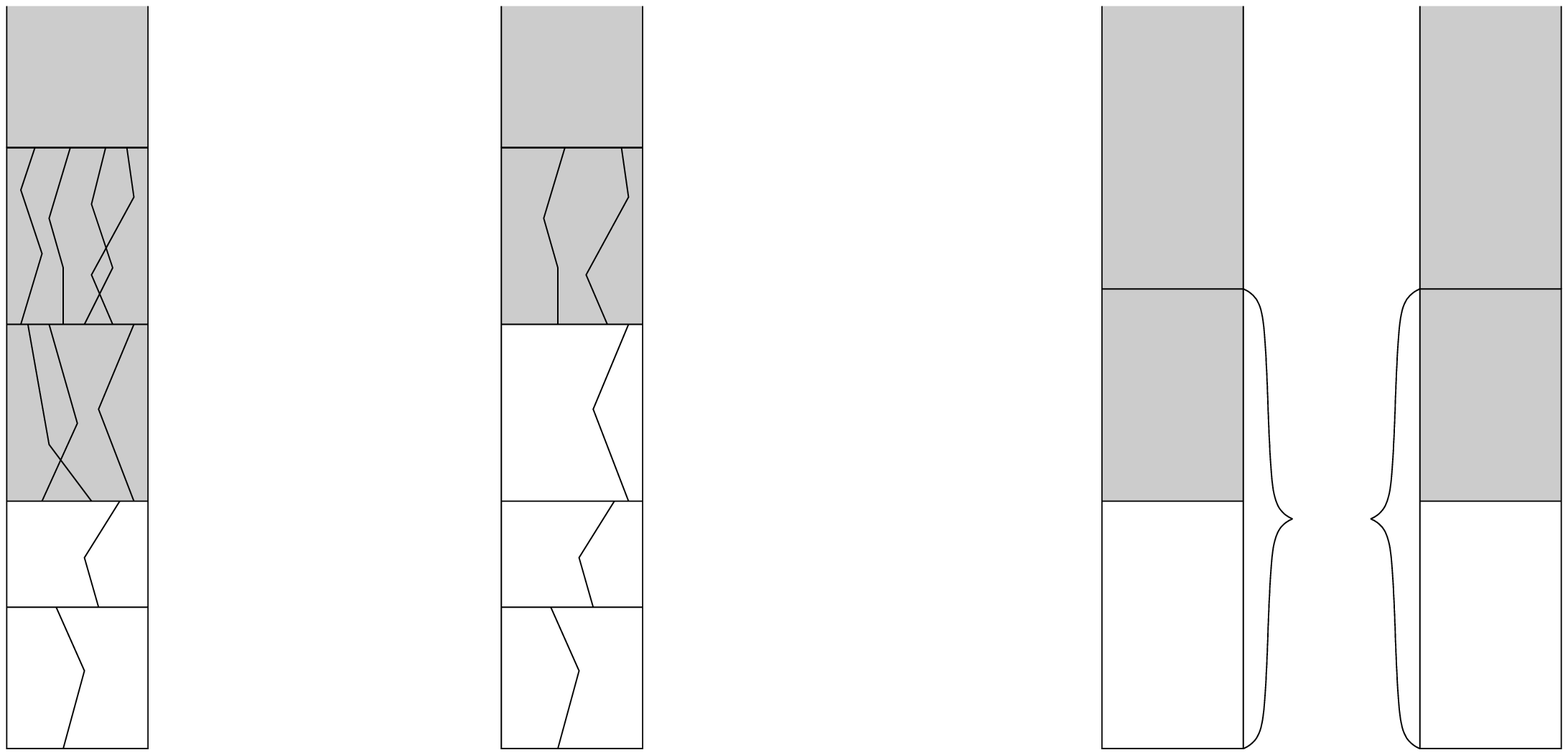}
    \caption{\label{fig:pq_new}  (a): $q\leq p$, $\trunklg(p)=2$, $\trunklg(q)=3$. (b): $q\leq_M p$}
  \end{center}
\end{figure}
So in particular $q\leq p$ implies that $\trunk(q)$ extends $\trunk(p)$, see
Figure~\ref{fig:pq_new}(a).

Of course we assume that there are sufficiently large creatures, otherwise
$\Qsu(\cK,\cS)$ is empty.\footnote{We need:
  For each $l\in\omega$ there is an $n\in\omega$ such that for all
  $m>n$ there is some $m$-creature with norm at least $l$.
}

The forcing $\Qsu(\cK,\cS)$ adds a generic real $\n\eta\DEFEQ \bigcup_{p\in
G}\trunk(p)$.  Note that when we have
halving (see next section), the generic filter $G$ is not
determined by $\n\eta$, at least not in the usual way.%
\footnote{If $\nor(\cc)$ is a function of $\val(\cc)$ and
  $\val(\cd)\subseteq\val(\cc)$ implies $\cd\in\cSc$, then the generic filter
  is determined by $\n\eta$. This assumption is reasonable
  (and is satisfied in many creature forcing constructions),
  but it is incompatible with halving.}

A note on the requirement
\begin{equation}\label{eq:jhi25}
  \nor(p(i))>0 \text{ for each }i\geq \trunklg(p)
\end{equation}
in the definition of $\Qsu$:
\begin{itemize}
  \item We could drop~\eqref{eq:jhi25}, since in the resulting 
    forcing notion the conditions that additionally satisfy~\eqref{eq:jhi25} are dense anyway.
  \item Because of~\eqref{eq:jhi25}, we are really only interested in creatures with
    norm $>0$, so we could restrict ourselves to creating pairs
    containing only such creatures.
  \item Alternatively, we could
    omit the concept of trunk from the definition altogether. Instead,
    we could assume the following: For all $\cc\in\cK(n)$
    and all $s\in\val(\cc)$ there is a $\cd\in \cS(\cc)$ such that
    $\val(\cd)=\{s\}$ (and therefore $\nor(\cd)=0$).
    However, this is not the ``right'' way to think about creature forcing,
    and this version could not be generalized to our
    variant of the countable support product.
\end{itemize}

In the rest of the section, we briefly comment on how our setting
fits into the framework of creature forcing developed
in~\cite{RoSh:470}:

A pair $(\cK,\cS)$ as above is a creating pair as defined
in~\cite[1.2]{RoSh:470}.  It satisfies the following additional properties:
\begin{itemize}
  \item finitary~\cite[1.1.3]{RoSh:470}: $\bfH(n)$ and  $\cS(\cc)$
    are always finite.
  \item simple~\cite[2.1.7]{RoSh:470}:
    $\cS$ is defined on single creatures
    only.\footnote{In non-simple creating pairs we can have something like
    $\cd\in \cS(\{\cc_1,\cc_2\})$, e.g.\
    $\cc_1$ could live on the interval $I_1$, $\cc_2$ on $I_2$,
    and $\cd$ is $\cc_1$ and $\cc_2$ ``glued together''.}
  \item forgetful~\cite[1.2.5]{RoSh:470}:
    $\val(\cc)$ does not depend on values of the
    generic real outside of the interval of $\cc$.\footnote{In
    the general case, $\val(\cc)$ is
    defined as a set of pairs $(u,v)$ where
    $v\in \prod_{i<\bfF (n+1)}\bfH(i)$ and \mbox{$u=v\restriction \bfF (n)$}.
    The intended meaning is
    that $\cc$ implies:
    If the generic object $\n\eta$ restricted to $\bfF (n)$ is $u$, then
    the possible values $v$ for $\n\eta\restriction \bfF (n+1)$ are
    those $v$ such that $(u,v)\in\val(\cc)$.
    Then ``$\cc$ is forgetful'' is defined as:
    If $(u,v)\in\val(\cc)$ and $u'\in \prod_{i<\bfF (n)}\bfH(i)$
    then $(u',v)\in\val(\cc)$. So in
    the forgetful
    case $\val(\cc)$ and $\{v:\, (\exists u)\,(u,v)\in\val(\cc)\}$
    carry the same information. In this paper
    we call the latter set $\val(\cc)$,  for simplicity of notation.}
  \item nice and smooth~\cite[1.2.5]{RoSh:470}:
    A technical requirement that is trivial in the case of forgetful
    simple creating pairs.
\end{itemize}

In~\cite{RoSh:470} two main frameworks for forcings are examined:
creature forcings~\cite[1.2.6]{RoSh:470} (defined by a creating pair~\cite[1.2.2]{RoSh:470})
and tree creature forcings~\cite[1.3.5]{RoSh:470}
(defined via a tree-creating pair~\cite[1.3.3]{RoSh:470}).
So in this paper we deal with creature forcings.\footnote{Actually
  every simple forgetful creating pair can be interpreted as
  tree-creating pair as well. The resulting
  tree-forcing however is different from the creature forcing:
  the creature forcing corresponds to the ``homogeneous'' trees only.}

In~\cite{RoSh:470} several ways to define forcings from a creating pair are
introduced. One example is lim-sup creature forcing $\Qswu$ defined
in~\cite[1.2.6]{RoSh:470}.  {\em Many simple cardinal invariants}~\cite{GoSh:448}
uses (a countable support product of) such forcings.  The lim-inf case $\Qsu$
is generally harder to handle, and~\cite[1.4.5]{RoSh:470} proves
that $\Qsu$ can collapse $\om1$.  In the rest of~\cite{RoSh:470},
$\Qsu$ is only considered in a special case (incompatible with simple) where
$\Qsu$ is actually equivalent to other forcings that are better behaved
(cf.~\cite[p23 and 2.1.3]{RoSh:470}).  We will introduce additional assumptions
(increasingly strong bigness and halving) to guarantee that $\Qsu$ is proper
and $\omega\ho$-bounding. These assumptions will actually make $\Qsu$
similar to $\Qsf$ of~\cite{RoSh:470}.

\section{bigness and halving, properness of $\Qsu$}\label{sec:puredecondim}

We will now introduce properties that guarantee that $\Qsu$ is proper.

\begin{Def}\label{def:big}
  Let $0<r\leq 1$, $B\in\omega$.
  \begin{itemize}
  \item $\cc$ is $(B,r)$-big
      if for all functions $F:\val(\cc)\fnto B$ there is a
      $\cd\in\cSc$
      such that $\nor(\cd)\geq \nor(\cc)-r$ and
      $F\restriction\val(\cd)$ is constant.\footnote{This
      is a variant of, but technically not quite the same as,
      \cite[2.2.1]{RoSh:470}.}
  \item $\cK(n)$ is $(B,r)$-big if every $\cc\in \cK(n)$
      with $\nor(\cc)>1$ is $(B,r)$-big.
  \item $\cc$ is $r$-halving,\footnote{cf.~\cite[2.2.7]{RoSh:470}.
    The original definition used
    $\nor(\chalf(\cc))\geq \nor(\cc)/2$ instead of $\nor(\cc)-r$,
    therefore the name halving.}
    if there is a $\chalf(\cc)\in\cSc$ such that
    \begin{itemize}
    \item $\nor(\chalf(\cc))\geq \nor(\cc)-r$, and
    \item if $\cd\in\cS(\chalf(\cc))$ and $\nor(\cd)>0$, then
      there is a $\cd'\in\cSc$ such that\\
      $\nor(\cd')\geq \nor(\cc)-r$ and $\val(\cd')\subseteq\val(\cd)$.
    \end{itemize}
  \item $\cK(n)$ is $r$-halving, if all $\cc\in\cK(n)$ with $\nor(\cc)>1$
    are $r$-halving.
  \end{itemize}
\end{Def}
So given $\cc$ and $\cd\in\cS(\chalf(\cc))$ as in the definition
of halving, we can ``un-halve'' $\cd$ to get $\cd'$.
Note that this $\cd'$ generally is not in $\cS(\chalf(\cc))$, although
$\val(\cd')\subseteq \val(\cd)\subseteq\val(\chalf(\cc))$.

Every creature is $(1,r)$-big.
If  $r'$ is smaller than $r$, then
$(B,r')$-bigness
implies $(B,r)$-bigness, and
$r'$-halving implies $r$-halving. We also get:
\begin{equation}\label{eq:valbig}
  \text{If }\cc\text{ is }(B,r)\text{-big and }0<r<\nor(\cc)\text{, then }
  B<|\val(\cc)|.
\end{equation}
An example for creatures with bigness and halving (and the much stronger
property decisiveness) can be found in Section~\ref{sec:example}.

We now show that increasing bigness and halving implies properness:
\begin{Thm}\label{thm:onedim}
  Set
  $\varphi(\mathord< n)\DEFEQ \prod_{i <\bfF (n)}\bfH(i)$
  and $r(n)\DEFEQ 1/(n\varphi(\mathord< n))$.
  If $\cK(n)$ is $(2,r(n))$-big and
  $r(n)$-halving for all $n$,
  then $\Qsu(\cK,\cS)$ is $\omega\ho$-bounding and proper
  and preserves the size of the continuum (in the following sense:
  in the extension, there
  is a bijection between the reals and old reals).
\end{Thm}
So in particular, CH is preserved.

\begin{Note}
  Only the growth rate of $r$ is relevant here. In particular: Fix some
  $\delta>1$. Then the theorem remains valid if we replace $(2,r(n))$-big
  and $r(n)$-halving with the weaker 
  condition $(2,\delta\cdot r(n))$-big and $\delta\cdot
  r(n)$-halving.  Also, it does not make any difference if we require bigness and
  halving only for those creatures with norm bigger than $\delta$ (instead of for
  all creatures with norm bigger than $1$).
\end{Note}

Note that $\varphi(\mathord< n)$ is the number of possible 
values for $\n\eta\restriction \bfF(n)$, or equivalently
the number of possible trunks with trunk-length $n$.

We also set  $\varphi(\mathord\leq n)=\varphi(\mathord< n+1)$ and
$\varphi(\mathord=n)=\varphi(\mathord\leq n)/\varphi(\mathord< n)=
\prod_{\bfF (n)\leq i <\bfF (n+1)}\bfH(i)$.

In the rest of this section we set $P=\Qsu(\cK,\cS)$.

We use a standard pure decision argument:

Let $\val(p,\mathord < n)$ denote $\Pi_{i<n}\val(p(i))$, the set of possible
values (modulo $p$) for $\n \eta\restriction \bfF(n)$.
The size of this set is at most $\varphi(\mathord<n)$.

We define for every $s\in \Pi_{i<\bfF(n)}\bfH(i)$
a condition $p\wedge s$:
$\trunklg(p\wedge s)=\max(n,\trunklg(p))$, and 
\[
  (p\wedge s)(i)\DEFEQ
    \begin{cases}
      s\restriction [F(i),F(i+1)-1]  &  \text{if }i<n\\
      p(i)  &  \text{otherwise.}
    \end{cases}
\]
We use this notion mostly for $s\in \val(p,\mathord < n)$. In this case,
$p\wedge s\leq p$. Note that 
\begin{equation}\label{eq:x}
  \{p\wedge s:\, s\in \val(p,\mathord < n)\}\text{ is predense under }p,
\end{equation}
which implies for all $s\in \val(p,\mathord < n)$
\begin{equation}\label{eq:gut}
  p\wedge s\, \forc\, \varphi\text{ iff }p\,\forc\, (s<\n\eta \rightarrow \varphi).
\end{equation}

$q\leq^* p$ means that $q$ forces $p$ to be in the generic filter.
\begin{equation}\label{eq:gorx}
  q\leq^* p\text{ implies }\val(q,\mathord < n)\subseteq  \val(p,\mathord < n).
\end{equation}
It is important to note that $\val(q(i))\subseteq \val(p(i))$ for all $i$ does
{\em not} imply $q\leq^* p$ (or even just $q\comp p$), since 
$\val(\cd)\subseteq\val(\cc)$ does not imply $\cd\in\cS(\cc)$.
(This would contradict halving.) However, the following {\em does}\/
follow from~\eqref{eq:x}:
\begin{equation}\label{eq:gur}
  \text{If $\val(q(i))\subseteq \val(p(i))$ for all $i\leq h$ and $q(i)\in
  \cS(p(i))$ for all $i> h$, then $q\leq^* p$.}
\end{equation}

Let $\n\tau$ be a name of an ordinal. $p$ $\mathord<n$-decides $\n\tau$, if
$p\wedge s$ decides\footnote{i.e.\ there is an $\alpha_s\in V$ such that
  $p\wedge s$ forces $\n\tau=\std \alpha_s$.}
$\n\tau$ for all $s\in\val(p,\mathord < n)$.  $q$ essentially decides $\n\tau$,
if $p$ $\mathord<n$-decides $\n\tau$ for some $n$.

So if $p$ essentially decides $\n\tau$, then we can calculate
the value of $\n\tau$ from a finite set of possible trunks of
$p$. So~\eqref{eq:gut} and~\eqref{eq:gorx} imply:
\begin{equation}\label{eq:y0}
  \text{If $p$ $\mathord<n$-decides $\n\tau$, and $q\leq^* p$,
  then $q$ $\mathord<n$-decides $\n\tau$.}
\end{equation}
We also get:
\begin{equation}\label{eq:y}
  \text{If $q\wedge s$ essentially decides $\n\tau$ for each
  $s\in\val(q,\mathord<n)$, then so does $q$.}
\end{equation}

We define the following (non-transitive) 
relations $\leq_n$ ($n\in\omega$) on $P$:
\begin{equation}
  \parbox{0.8\columnwidth}{
    \raggedright
    $q\leq_n p$ if
    $q\leq p$ and there is an $h\geq n$ such that
    $q\restriction h=p\restriction h$
    and $\nor(q(i))\geq n$ for all $i\geq h$.
  }
\end{equation}
(Cf.~Figure~\ref{fig:pq_new}(b) on page \pageref{fig:pq_new}).

\begin{proof}[Proof of Theorem~\ref{thm:onedim}]
We will show the following properties:
\begin{itemize}
  \item $q\leq_0 p$ implies $q \leq p$, and $q\leq_{n+1} p$ implies $q\leq_n p$.
  \item {\em (Fusion.)}\/ 
    For every sequence $p_0\geq_0 p_1\geq_1 p_2\geq \dots$
    there is a $q$ stronger than each $p_n$.
  \item {\em (Pure decision.)}\/
     For every name $\n \tau$ of an ordinal, $n\in\omega$,
    and $p\in P$, there is a $q\leq_n p$ essentially deciding $\n\tau$.
\end{itemize}

Then the standard argument can be employed to 
show Theorem~\ref{thm:onedim}:
\begin{itemize}
  \item {\em $\omega\ho$-bounding:}
    Let $\n f$ be the name for
    a function from $\omega$ into ordinals and $p\in P$.
    Set $p_0=p$.
    If $p_n$ is already constructed, choose $p_{n+1}\leq_{n+1}p_n$
    essentially deciding $\n f(n)$. 
    Fuse the sequence into some $q$. Then modulo $q$ there
    are only finitely many possibilities for each $\n f(n)$.
  \item
  {\em Proper:}
    Let $N\esm H(\chi)$ be countable and contain $P$ and $p_0$.
    Let
    $(\n \tau_n)_{n\in\omega}$ list the
    $P$-names of ordinals that are in $N$.
    Choose (in $N$) $p_{n+1}\leq_n p_n$ such that $p_{n+1}$
    essentially decides $\n\tau_n$. If $q\leq p_n$ for all $n$,
    then $q$ is is $N$-generic.
  \item
  {\em The size of the continuum:}
    So for every $p$ in $P$ and $P$-name $\n r$ for a real
    there is a $q\leq p$ continuously reading $\n r$.
    This means that $\n r$ is calculated by a function 
    \[
      \text{eval}:\bigcup_{n\in\omega}\val(q,\mathord<n)\to 2^{<\omega}.
    \]
    (Since each $\n r(m)$ is determined by $\val(q,\mathord<M)$
    for some $M$.)
    There are only $2^{\al0}$ many such functions, and $|P|=2^{\al0}$ many
    conditions.
\end{itemize}

So we just have to show pure decision and fusion. Fusion is easy:
Let $(p_n)_{n\in\omega}$
satisfy $p_{n+1}\leq_{n+1} p_n$.
Set $q(n)=p_{n}(n)$. Then
$q$ is in $P$: Fix any $M\in\omega$. There
is an $h>M$ such that
\begin{equation}\label{eq:few}
  \nor(p_M(m))\geq M\text{ for all }m\geq h.
\end{equation}
Then~\eqref{eq:few} holds for $p_{M+1}$ as well,
and for each $p_k$ with $k>M$, and therefore for $q$. Clearly,
$q\leq p_n$ for each $n$.

It remains to be shown that $P$ satisfies
pure decision.

Let $\n\tau$ be the name of an ordinal.

{\bf The basic construction $S(p,M)$:}\\
  Assume that $\trunklg(p)=n$ and $M\in\omega$.
  We define $S(p,M)$ the following way, see Figure~\ref{fig:pbasic}:

\begin{figure}[tb]
  \begin{center}
    \newcommand{\mydecr}{$\mathord\geq \cdot-r$}
    \scalebox{0.4}{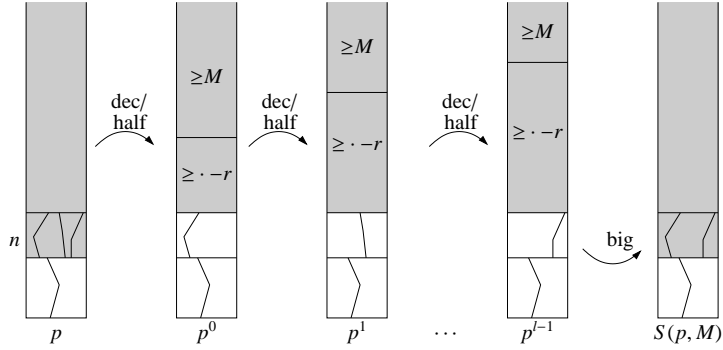}
    \caption{\label{fig:pbasic} 
    The basic construction $S(p,M)$.}
  \end{center}
\end{figure}

  Enumerate $\val(p,\mathord\leq n)$ as $s^0,\dots,s^{l-1}$.
  So $l\leq \varphi(\mathord=n)$.
  Set $p^{-1}=p$. Given $p^k$, define $p^{k+1}\in P$ as follows:
  $\trunk(p^{k+1})=s^{k+1}$, $p^{k+1}\leq p^k\wedge s^{k+1}$,
  and there is an $h^{k+1}$ such that
  \begin{itemize}
    \item if $n<m<h^{k+1}$, then $\nor(p^{k+1}(m))\geq \nor(p^{k}(m))-r(m)$,
    \item if $m\geq h^{k+1}$, then $\nor(p^{k+1}(m))\geq M$,
  \end{itemize}
  and such that additionally one of the following two cases holds:%
  \begin{description}
    \item[dec]  $p^{k+1}$ essentially decides $\n\tau$, or
    \item[half] it is not possible to satisfy ``dec'' (for any
      choice of $h^{k+1}$),
      then $p^{k+1}(m)=\chalf(p^{k}(m))$ for all $m>n$.
  \end{description}
  This way we construct $p^k$ for each $0\leq k<l$.
  At each step $0\leq k<l$, we have one of the cases ``dec'' or ``half''.
  This gives a function $F: \val(p(n))\to \{\text{dec},\text{half}\}$,
  and we use bigness to thin out $p(n)$ and
  get some $\cd\in\cS(p(n))$ such that $F\restriction\val(\cd)$ is constant
  and $\nor(\cd)\geq \nor(p(n))-r(n)$.

  Note that in this construction we have to assume that 
  $\nor(p^{k}(m))>1$ for all $-1\leq k<l-1$ and $m>n$, otherwise we cannot 
  halve $p^{k}(m)$. Also, 
  $\nor(p(n))$ has to be bigger than $1$, otherwise we
  cannot use bigness. 
  Let $S(p,M)$ be undefined if these conditions are not met.
  Otherwise,
  we define $q=S(p,M)$ as follows: 
  \[
    q\restriction n=p\restriction n=\trunk(p),\quad
    q(n)=\cd,\quad q(m)=p^{l-1}(m)\text{ for }m>n.
  \]

  We call $q$ halving, if the constant value of $F\restriction \val(q(n))$ 
  is ``half''.
  We will show that $q$ cannot be halving. 

  If $q$ is not halving, i.e.\ if the constant value is ``dec'', then
  $q$ essentially decides $\n\tau$: 
  If $s\in \val(q,\mathord\leq n)$, then $s=s^k$ for some $k<l$, and
  $q\wedge s\leq p^k$ essentially decides $\n\tau$.
  Now use~\eqref{eq:y}.

  {\bf Some properties of  $S(p,M)$:}\\
  If $q=S(p,M)$ is defined, then it satisfies the following:
  \begin{gather}
    \label{eq:qa} \nor(q(n))\geq \nor(p(n))-r(n).\\
    \label{eq:qb} \text{If $m>n$, then }\nor(q(m))\geq
      \min(M,\nor(p(m)))-\varphi(\mathord= n)\cdot r(m).\\
    \label{eq:dreizwei}\text{If $q$ is halving, then no $q'\leq q$ with
      trunk-length $n+1$ essentially decides $\n\tau$.}
  \end{gather}
  To see~\eqref{eq:dreizwei}, assume that $q'$
  is a counterexample.
  So $q'\leq q\wedge s^k\leq p^{k}$ for some 
  $0\leq k<l$, and $\nor(q'(m))>0$ 
  for all $m>n$. Since $q$ is halving,
  $p^{k}$ was produced by halving $p^{k-1}$. Pick an $h$ such that 
  $\nor(q'(m))>M$ for all $m\geq h$. For $n<m<h$,
  $p^{k}(m)=\chalf(p^{k-1}(m))$ and $q'(m)\in \cS(q(m))\subseteq
  \cS(p^{k}(m))$,
  so we can un-halve 
  $q'(m)$ to get some $\cd_m\in \cS(p^{k-1}(m))$
  with $\val(\cd_m)\subseteq \val(q'(m))$ and 
  $\nor(\cd_m)\geq \nor(p^{k-1}(m))-r(m)$.
  But then we could have chosen a deciding condition $r$ instead
  of $p^{k}$: Define $r(m)=\delta_m$ for $n<m<h$ and
  $r(m)=q'(m)$ otherwise.
  According to~\eqref{eq:gur},
  \mbox{$r\leq^* q$. \eqref{eq:y0}}
  implies that $r$ essentially decides
  $\n\tau$, a contradiction.

  {\bf $S(p,M)$ essentially decides:}\\
  We assume that $S(p,M)$ is halving and get a 
  contradiction the following way:
  We show that the ``successors'' of $q$ with
  increased stem have to be halving as well, and we 
  can fuse them into some $q^\omega$. But there will
  be a $q'\leq q^\omega$ deciding $\n\tau$, a contradiction.
  In more detail:
  \begin{equation}\label{eq:xy7}
    \parbox{0.8\columnwidth}{
      \raggedright If $\trunklg(p)=n$, $\nor(p(m))>3$ for all $m\geq n$ and if $M>3$,
        then $S(p,M)$ exists and is not halving.}
  \end{equation}
  Assume towards a contradiction that
  $S(p,M)$ is halving (or does not exist).
  Set $q^{n-1}=p$.
  Assume that for $k\geq n-1$, we have
  already defined  $q^k$. We set
  $M_k=M+k+1-n$ (note that $M_{n-1}=M$), and define $q^{k+1}$ the 
  following way:


  List $\val(q^k,\mathord\leq k)$ as $s^0,\dots,s^{l-1}$. So
  $l\leq \varphi(\mathord \leq k)$.
  Set $r^{-1}=q^k$. Given $r^{i-1}$, 
  set 
  \begin{equation}\label{eq:rx}
    r^{i}=S(r^{i-1}\wedge s^i,M_k)
  \end{equation}
  (if defined). So $r^i$ has trunk-length $k+1$.
  Define $q^{k+1}(m)$ to be $q^{k}(m)$ for $m\leq k$ and
  $r^{l-1}(m)$ otherwise.

  So in particular, $q^n=S(p,M)$.

  If $q^{k+1}$ is defined, then~\eqref{eq:qa} and~\eqref{eq:qb} imply:
  \begin{itemize}
    \item $q^{k+1}(m)=q^{k}(m)$ for $m\leq k$.
    \item $\nor(q^{k+1}(k+1))\geq \nor(q^{k}(k))-
                              \varphi(\mathord\leq k)\cdot r(k+1)$.
    \item $\nor(q^{k+1}(m))\geq \min(M_k,\nor(q^k(m)))-
                              \varphi(\mathord\leq k+1)\cdot r(m)$
      for $m>k+1$.
  \end{itemize}
  So in any case, we get for all $m\in\omega$
  \begin{align}
    \label{eq:qc}
      \nor(q^{k+1}(m))\geq&\min(M_k,\nor(q^k(m)))-
                         \varphi(\mathord<m)\cdot r(m).
    \\\intertext{Iterating this $l$ many steps
      (note that $q^k(m)$ remains constant if $k\geq m$) we get for all $m$:}
    \label{eq:qf} 
      \nor(q^{k+l}(m))\geq & \min(M_k,\nor(q^k(m)))-\min(l,m-k)
      \cdot  \varphi(\mathord<m)\cdot r(m),
    \\\intertext{and since $r(m)=1/(m\cdot \varphi(\mathord<m))$, we get}
    \label{eq:qg} 
      \nor(q^{k+l}(m))\geq & \min(M_k,\nor(q^k(m)))-1.
  \end{align}
  If we set $k=n-1$, this shows that
  $\nor(q^{k+l}(m))\geq 2$ for all $l\in\omega$,
  and that therefore $q^{k+l+1}$ is defined.
  Also, if we 
  define $q^\omega$ by $q^\omega(m)=q^m(m)$, then
  $q^\omega\in P$:
  Given $N\in\omega$, pick 
  $k$ such that $M_k>N+1$ and pick
  $h>k$ such that $\nor(q^k(m))>N+1$
  for all $m>h$. If $m>h$, i.e.\ $m=k+l$ for some $l>0$,
  then $q^\omega(m)=q^{k+l}(m)$,
  and $\nor(q^{k+l}(m))\geq \min(M_k, \nor(q^k(m))-1>N$.

  Also, $q^\omega\leq q^k$ for all $k\in\omega$.

  The property~\eqref{eq:dreizwei} of $S$ 
  can by induction be generalized to any $k\geq n$
  (recall that $q=S(p,M)=q^n$).
  \begin{equation}\label{eq:prop7}
    \text{No $q'\leq q^{k}$ with trunk-length $k+1$ essentially
      decides $\n\tau$.}
  \end{equation}
  For $k=n$ this is~\eqref{eq:dreizwei}. 
  We assume that~\eqref{eq:prop7} holds for $k$ and show it 
  for $k+1$. Assume $q'$ is a counterexample. $q'$ is
  stronger than some of the $r^i$ ($0\leq i<l$) used 
  in~\eqref{eq:rx} to construct $q^{k+1}$. 
  $r^i=S(r^{i-1}\wedge s^i,M_k)$ has trunk-length $k+1$
  and is stronger than $q^k$, so
  we can apply~\eqref{eq:prop7} to see that $r^i$ cannot 
  essentially decide $\n\tau$. So $r^i$ is halving. 
  Using~\eqref{eq:dreizwei}, we see that no $q'\leq r^i$
  with trunk-length $k+2$ essentially decides $\n\tau$,
  a contradiction.

  On the other hand, there is a $q'\leq q^\omega$ deciding $\n\tau$.
  Set $k=\trunklg(q)-1$. Then $q'\leq q^\omega\leq q^k$
  contradicts~\eqref{eq:prop7}.

  {\bf Pure decision:}\\
  Given $p\in P$ and $M\in\omega$, pick
  $n$ such that $p(m)>M+5$ for all $m\geq n$.
  Similarly to above, enumerate $\val(p,\mathord<n)$
  as $s^0,\dots,s^{l-1}$, set $r^{-1}=p$ and $r^{k+1}=S(r^k\wedge s^{k+1},M+5)$.
  Define $q$ by $q\restriction n=p\restriction n$ and
  $q(m)=r^{l-1}(m)$ for $m\geq n$.
  Just as in~\eqref{eq:qc},
  $\nor(q(m))\geq \min(M+5,\nor(p(m)))-1> M+4$ for $m>n$, i.e.\ $q\leq_M p$.
  As we already know by~\eqref{eq:xy7}, each $r^k$ essentially
  decides $\n\tau$, so by~\eqref{eq:y},
  $q$ essentially decides $\n\tau$ as well.
\end{proof}

A simple modification of the proof leads to a stronger property:
Using the same $\varphi$ and $r$ as in the previous theorem, we get:
\begin{Thm}
  Assume that $g:\omega\to\omega\setminus 1$ is 
  monotonously increasing, that $\n\nu$ is a $P$-name 
  and that $p\in P$ forces that $\n\nu(n)<g(n)$ for all $n$.
  If each $\cK(n)$ is $(g(n)+1,r(n))$-big and $r(n)$-halving,
  then there is a $q\leq p$ which 
  $\mathord<n$-decides $\n\nu(n)$ for all $n$.
\end{Thm}
We call this phenomenon {\em rapid reading}.

\begin{proof}
  We modify the last proof in the following way:

  {\em The basic construction $S(p,l,M)$:}
  We again assume that $n=\trunklg(p)$, and use the notation
  $S(p,l,M)$ (for $l\leq n$)
  for the same construction as $S(p,M)$, where 
  we set $\n\tau=\n\nu(l)$, and
  instead of trying to {\em essentially} decide $\n\tau$,
  we try to {\em decide} it.
  So 
  instead of the two cases ``dec'' and ``half'', we get 
  $g(l)+1$ many cases: ``0'', \dots, ``$g(l)-1$'', and
  (if none of these cases can be satisfied) ``half''.
  Since $l\leq n$ and $g$ is increasing, we can use
  $(g(n)+1,r(n))$-bigness instead of just
  $(2,r(n))$-bigness, and we again get a homogeneous $\cd$.
  If $S(p,l,M)$ is not halving, then it decides $\n\nu(l)$.

  {\em Some properties of  $S(p,l,M)$:}
  We again get~\eqref{eq:qa} and~\eqref{eq:qb}, and in~\eqref{eq:dreizwei}
  we replace ``essentially decides $\n\tau$'' with
  ``decides $\n\nu(l)$'', i.e.\ we get:
  \begin{equation*}
    \text{If $q$ is halving,
    then no $q'\leq q$ with trunk-length $n+1$ decides $\n\nu(l)$.}
  \end{equation*}

  {\em $S(p,l,M)$ decides:} We again construct $q^k$, each time trying to
  decide $\n\tau=g(l)$ (independently of $k$). 
  So~\eqref{eq:rx} now reads:
  \begin{equation*}
    r^i=S(r^{i-1}\wedge s^i,l,M_k).
  \end{equation*}
  (Here we only need $(g(l)+1,r(k))$-bigness). 
  Again we get~\eqref{eq:qg}, and  therefore
  each $q^k$ (and $q^\omega$) is
  defined, and~\eqref{eq:prop7} now tells us
  \begin{equation*}
    \text{No $q'\leq q^{k}$ with trunk-length $k+1$ decides 
  $\n\tau$.}
  \end{equation*}
  But there is some $q'\leq q^\omega$ deciding $\n\tau$,
  a contradiction.

  So far we know the following: 
  \begin{equation}\label{eq:gube7}
    \parbox{0.8\columnwidth}{
      \raggedright If $\trunklg(p)=n$, $\nor(p(m))>3$ for $m\geq n$, and $M>3$,
      then $S(p,n,M)$ exists and decides $\n\nu(n)$.}
  \end{equation}
  {\em Rapid reading:}
  Instead of the part on {\em pure decision}, we proceed as follows:
  Given $p\in P$, we can assume (by enlarging the stem)
  that 
  $\nor(p(m))>5$ for all $m>\trunklg(p)$.
  We set $k_0=\trunklg(p)-1$ and
  $q^{k_0}=p'$.
  We now construct $q^k$ and $q^\omega$ just as above, but this 
  time using 
  \begin{equation*}
    r^{i}=S(r^{i-1}\wedge s^i,k+1,M_k).
  \end{equation*} 
  As in~\eqref{eq:qg} we see that $r^i$, $q^k$ and $q^\omega$ exist.
  $r^i$ has sufficient norm and trunk-length $k+1$, so by~\eqref{eq:gube7}
  each $r^i$ 
  decides $\n \nu(k+1)$. This implies that $q^{k+1}$ (and therefore
  $q^\omega$ as well) $\mathord\leq k$-decides $\n \nu(k+1)$.
\end{proof}

Note that $P$ has size continuum, and in particular it is
$(2^{\al0})^+$-cc. Together with proper, that gives us:
\begin{Lem}
  Under CH and the assumptions of
  Theorem~\ref{thm:onedim}, $P$ preserves all cardinals (and cofinalities)
  and the size of the continuum.
\end{Lem}

\section{decisiveness, properness of finite products}\label{sec:decisive}

In this section, we fix a {\em finite\/} set $I$
and for every $i\in I$ a
creating pair $(\cK_i,\cS_i)$.

The product forcing $\prod_{i\in I} \Qsu(\cK_i,\cS_i)$ is equivalent to
$\Qsu(\cK_I,\cS_I)$, where the creating pair $(\cK_I,\cS_I)$
is defined as follows:
An $n$-creature $\cc\in\cK_I(n)$
corresponds to an $|I|$-tuple $(\cc_i)_{i\in I}$ such that
$\cc_i\in \cK_i(n)$.
$\val(\cc)=\prod_{i\in I}\val(\cc_i)$,
$\nor(\cc)=\min(\{\nor(\cc_i):\, i\in I\})$, and
$\cd=(\cd_i)_{i\in I}$ is in $\cS(\cc)$ if
$\cd_i\in\cS(\cc_i)$ for all $i\in I$.\footnote{So
an $n$-creature ``lives'' on the product $\prod_{i\in I}[\bfF _i(n),\bfF _i(n+1)-1]$. This does not fit our restrictive framework, so we could
just ``linearize'' the product. Assume $I\in\omega$, i.e.\ $I=\{0,\dots,I-1\}$.
Set
$\bfF _I(n)\DEFEQ\sum_{i\in I} \bfF _i(n)$ and write it in the following way:\\
\setlength{\unitlength}{5mm}
\newlength{\mul}\setlength{\mul}{0.9\unitlength}
\centerline{\begin{picture}(23,2)(-1,-1.2)
\put( 0,0){\line(1,0){21}}
\put( 0,-0.5){\line(0,1){1.0}}
\put(0,-1){\makebox(0,0){\footnotesize $\bfF _I(0)$}}
\put( 3,-0.15){\line(0,1){0.3}}
\put( 1.5,-0.7){\makebox(0,0){\footnotesize $\underbrace{\hspace{3\mul}}_{\bfF _0(1)}$}}
\put( 5,-0.15){\line(0,1){0.3}}
\put( 4,-0.7){\makebox(0,0){\footnotesize $\underbrace{\hspace{2\mul}}_{\bfF _1(1)}$}}
\put( 7,-0.5){\makebox(0,0){\footnotesize \dots}}
\put( 9,-0.15){\line(0,1){0.3}}
\put( 11,-0.7){\makebox(0,0){\footnotesize $\underbrace{\hspace{4\mul}}_{\bfF _{I-1}(1)}$}}
\put( 13,-0.5){\line(0,1){1.0}}
\put( 13,-1){\makebox(0,0){\footnotesize $\bfF _I(1)$}}
\put( 16,-0.15){\line(0,1){0.3}}
\put( 14.5,-0.7){\makebox(0,0){\footnotesize $\underbrace{\hspace{3\mul}}_{\bfF _0(2)}$}}
\put( 18,-0.5){\makebox(0,0){\footnotesize \dots}}
\end{picture}}
Now it should be clear how to formally define $\bfH_I$, $\cK_I$, $\cS_I$ etc.}

If each $\cK_i(n)$ is $r$-halving, then $\cK_I(n)$ is $r$-halving as well: 
We can set $\chalf(\cc)\DEFEQ (\chalf(\cc_i))_{i\in I}$.  This satisfies 
Definition~\ref{def:big} of halving: Assume that  $\cd\in\cS( \chalf(\cc))$ and
$\nor(\cd)>0$. So $\cd=(\cd_i)_{i\in I}$, $\cd_i\in \cS(\cc_i)$, and
$\nor(\cd_i)>0$ for all $i\in I$.  We can un-halve each $\cd_i$
to some $\cd'_i$, and set $\cd'=(\cd'_i)_{i\in I}$. Then
$\cd'$ is as required.

However, $\cK_I$ will not satisfy bigness, since  a function
$F:\prod_{i\in I}\val(\cc_i)\to 2$ can generally not be
written as a product of functions $F_i:\val(\cc_i)\to 2$.
So to handle bigness we have to introduce a new notion:

\begin{Def} Let $0<r\leq 1$, $B,K,n>0$.
  \begin{itemize}
  \item $\cc$ is hereditarily $(B,r)$-big, if every $\cd\in \cSc$
    with $\nor(\cd)>1$ is $(B,r)$-big.
  \item $\cc$ is $(K,n,r)$-decisive, if
    there are $\cd^-,\cd^+\in \cSc$
    such that\\ $\nor(\cd^-),\nor(\cd^+)\geq \nor(\cc)-r$,
    $\card{\val(\cd^-)}\leq K$ and $\cd^+$ is hereditarily
    $(2^{K^n},r)$-big.\\
    $\cd^-$ is called a $K${\em -small successor}, and
    $\cd^+$ a $K${\em-big successor} of $\cc$.
  \item $\cc$ is $(n,r)$-decisive if $\cc$ is $(K,n,r)$-decisive for some $K$.
  \item $\cK(n)$ is  $(n,r)$-decisive if every $\cc\in \cK(n)$ with
    $\nor(\cc)>1$ is  $(n,r)$-decisive.
  \end{itemize}
\end{Def}

An example for decisive, halving creatures can be found in Section~\ref{sec:example}.

\begin{Lem}\label{lem:increasebigness}
  \begin{enumerate}
    \item
      If $\cc$ is $(n,r)$-decisive
      (i.e.\ $\cc$ is $(K_0,n,r)$-decisive for {\bf some} $K_0$),
      then for {\bf every} $K\in\omega$ there is 
      {\bf either} a $K$-big successor
      {\bf or} a $K$-small successor of $\cc$.
    \item
      If $\cc$ is $(K,n,r)$-decisive and hereditarily $(B,r)$-big,
      and if $\nor(\cc)>1+r$, then $B<K$.
    \item
      Assume that $\cK(n)$ is $(n,r)$-decisive and
      $(B,r)$-big for some $B\geq 1$, that $\delta\in\omega$
      and that $\nor(\cc)>1+\delta\cdot r$.
      Then there is a hereditarily $(\EXP(B,n,\delta),r)$-big $\cd\in\cSc$ such
      that $\nor(\cd)\geq \nor(\cc)-\delta\cdot r$, where
      $\EXP(B,n,0)=B$ and $\EXP(B,n,m+1)=2^{\EXP(B,n,m)^n}$.
    \item In particular, if $\cK(n)$ is $(n,r)$-decisive and
      $\nor(\cc)>1+r$, then there is a hereditarily $(2,r)$-big
      $\cd\in\cS(\cc)$ such that $\nor(\cd)\geq \nor(\cc)-r$.
    \item We can avoid small sets without decreasing the
      norm too much: Assume that $\cK(n)$ is $(n,r)$-decisive and
      $(B,r)$-big for some $B\geq 1$, that $\delta\in\omega$
      and that $\nor(\cc)>1+(\delta+1)\cdot r$.
      If $X\subseteq \val(\cc)$ has size less than $\EXP(B,n,\delta)$,
      then there is a $\cd\in\cSc$ such that 
      $\nor(\cd)\geq \nor(\cc)-(\delta+1)\cdot r$ and
      $\val(\cd)$ is disjoint to $X$.
  \end{enumerate}
\end{Lem}

\begin{proof}
  (1): If $K\leq K_0$, use $\cd^-$, otherwise use $\cd^+$.
  (2): The $K$-small successor $\cd^-$ is $B$-big, and $|\val(\cd^-)|<K$.
  Now use~\eqref{eq:valbig}.
  (3): Set $\cd^+_0=\cc$. Assume that $\cd^+_i$
  is defined and has norm bigger than 1. So
  $\cd^+_i$ is decisive, i.e.\ there is a $K_i$ and
  a $K_i$-small successor $\cd^-_{i+1}$ and
  a $K_i$-big successor $\cd^+_{i+1}$. According to
  (2), $K_0>B$, and $K_{i+1}>2^{K_i^n}\geq \EXP(B,n,i+1)$.
  In particular, 
  $\cd^+_\delta$ is hereditarily $\EXP(B,n,\delta)$-big. 
  (4): Every creature is $(1,r)$-big.
  (5) follows from (3): First get a $(\EXP(B,n,\delta),r)$-big
  creature $\cd_0$, then use the function $F$ that maps 
  $\val{\cd_0}$ to $X\cup\{\text{\tt NotInX}\}$ and thin
  out $\cd_0$ to get an $F$-homogeneous $\cd$.
\end{proof}

We now show by induction on $k$: If the $n$-creatures are $(k,r)$-decisive, then we can
generalize bigness to $k$-tuples.
\begin{Lem}\label{lem:decisive}  Assume that $k,m,t\geq 1$, $0<r\leq 1$,
  $\cc_0,\dots,\cc_{k-1}\in \cK(n)$ and $F$ satisfy the following:
  \begin{itemize}
  \item $\nor(\cc_i)>1+r\cdot (k-1)$,
  \item $\cK(n)$ is $(k,r)$-decisive and each $\cc_i$ is hereditarily
    $(2^{m^t},r)$-big, and
  \item $F$ is a function from $\prod_{i\in k} \val(\cc_i)$ to $2^{m^t}$.
  \end{itemize}
  Then there are $\cd_0,\dots,\cd_{k-1}\in \cK$ such that
  \begin{itemize}
  \item $\cd_i\in\cS(\cc_i)$,
  \item $\nor(\cd_i)\geq \nor(\cc_i)-r\cdot k$, and
  \item $F\restriction \prod_{i\in k} \val(\cd_i)$ is constant.
  \end{itemize}
\end{Lem}

\begin{proof}
  The case $k=1$ follows directly from Definition~\ref{def:big} of
  $(2^{m^t},r)$-big (decisive is not needed).
  So assume the lemma holds for $k$, and let us investigate the case $k+1$.

  $\cc_k$ is $(k+1,r)$-decisive, i.e.\ there is an $M$ such that
  $\cc_k$ is $(M,k+1,r)$-decisive. So~\ref{lem:increasebigness}(2) implies
  \begin{equation}\label{eq:Mx}
    M>2^{m^t}.
  \end{equation}
  According to~\ref{lem:increasebigness}(1), for 
  each $\cc_i$ ($i<k$) we can pick some $\cd_i$ that is either 
  an $M$-small successor or an $M$-big successor of $\cc_i$ (since each
  $\cc_i$ is $(k+1,r)$-decisive).
  If $\cd_0$ is $M$-small, then we let
  $\cd_k$ be the $M$-big successor of $\cc_k$, otherwise
  the $M$-small one. (For $\cc_k$ we have both options, since 
  $\cc_k$ is $(M,k+1,r)$-decisive.)

  This gives a sequence $(\cd_i)_{i\in k+1}$ satisfying $\cd_i\in
  \cS(\cc_i)$ and $\nor(\cd_i)\geq \nor(\cc_i)-r$.
  Set $S\DEFEQ\{i\in k+1:\, \cd_i \text{ is }M\text{-small}\}$, and
  $L\DEFEQ (k+1)\setminus S$.  
  So $\{L,S\}$ is a non-trivial partition of $k+1$, since
  $0$ and $k$ are in different sets.
  If $i\in S$, then $|\val(\cd_i)|<M$, if $i\in L$ then
  $\cd_i$  is hereditarily $2^{M^{k+1}}$-big.

  Set $Y\DEFEQ \prod_{i\in S}\val(\cd_i)$. $\card{Y}\leq M^{\card{S}}$.
  So we
  can write $Y$ as $\{y_1,\dots ,y_{M^{\card{S}}}\}$.

  Define $F^*$ on $\prod_{i\in L}\val(\cd_i)$ by
  \[
    F^*(x)\DEFEQ (F(x^\frown y_1),\dots,
    F(x^\frown y_{M^{\card{S}}})).
  \]
  So (using~\eqref{eq:Mx} for the last inequality) we get:
  \[
    \card{\image(F^*)}\leq
    \card{\image(F)}^{M^{\card{S}}}\leq 2^{m^t M^{\card{S}}}<
    2^{M^{\card{S}+1}}.
  \]

  For $i\in L$, $\cd_i$ is hereditarily $2^{M^{k+1}}$-big
  and therefore $2^{M^{|S|+1}}$-big, and $|L|\leq k$.
  Therefore we can apply the induction hypothesis to $k'\DEFEQ \card{L}$,
  $m'\DEFEQ M$, $t'\DEFEQ \card{S}+1$,
  $F'\DEFEQ F^*$ and $\cc'_i\DEFEQ \cd_i$ for $i\in L$.
  This gives us $(\cd'_i)_{i\in L}$ such that
  \begin{itemize}
  \item $\cd'_i\in\cS(\cd_i)\subseteq \cS(\cc_i)$,
  \item $\nor(\cd'_i)\geq
    \nor(\cd_i)-r\cdot k'\geq \nor(\cc_i)-r(k+1)$, and
  \item $F^*\restriction \prod_{i\in L}\val(\cd'_i)$ is constant, say
    $(F^{**}(y_1),\dots,F^{**}(y_{M^{\card{S}}}))$.
  \end{itemize}
  $F^{**}$ is a function from $Y=\prod_{i\in S}\val(\cd_i)$ to $2^{m^t}$.
  Now we apply the induction hypothesis again, this time to
  $k''\DEFEQ \card{S}<k+1$, $m''\DEFEQ m$,
  $t''=t$, $F''\DEFEQ F^{**}$, and $\cc''_i\DEFEQ \cd_i$ for $i\in S$.
  This gives us $(\cd'_i)_{i\in S}$ such that
  \begin{itemize}
  \item $\cd'_i\in\cS(\cd_i)\subseteq \cS(\cc_i)$,
  \item $\nor(\cd'_i)\geq \nor(\cd_i)-r\cdot k''\geq \nor(\cc_i)-r(k+1)$, and
  \item $F^{**}\restriction \prod_{i\in S}\val(\cd'_i)$ is constant.
  \end{itemize}
  Then $(\cd'_i)_{i\leq k}$ is as required.
\end{proof}

According to~\ref{lem:increasebigness}(3), we
can increase the hereditary bigness by decreasing the norm.
So we get (again setting 
      $\EXP(B,n,0)=B$ and $\EXP(B,n,m+1)=2^{\EXP(B,n,m)^n}$):
\begin{Cor}\label{cor:biggerbigness}
  Fix $\delta\geq 1$.
  Assume that $k\geq 1$, $0<r\leq 1$,
  $\cK(n)$ is $(k,r)$-decisive and $(B,r)$-big,
  $\nor(\cc_i)>1+r\cdot (\delta+k-1)$ for $0\leq i < k$
  and  $F: \prod_{i\in k} \val(\cc_i)\to \EXP(B,k,\delta)$.
  Then there are $\cd_i\in \cS(\cc_i)$ with $F$-homogeneous product
  such that $\nor(\cd_i)\geq \nor(\cc_i)-r\cdot (\delta+k)$.
\end{Cor}

\begin{proof}
  By first decreasing the norms by at most $\delta\cdot r$,
  we can assume that each $\cc_i$ is hereditary 
  $\EXP(B,k,\delta)$-big.
  Now use Lemma~\ref{lem:decisive}.
  (Note that $\EXP(B,n,\delta)$ is of the form $2^{m^t}$ for some $m$ and $t$.) 
\end{proof}

Every creature is $(1,r)$-big, and $\EXP(1,n,1)=2$. So we 
get for $\delta=1$:
\begin{Cor}\label{cor:multidimbignesssimple}
  Assume that $k\geq 1$, $0<r\leq 1$,
  $\cK(n)$ is $(k,r)$-decisive,
  $\nor(\cc_i)>1+r\cdot k$ for $0\leq i < k$
  and $F: \prod_{i\in k} \val(\cc_i)\to 2$.
  Then there are $F$-homogeneous $\cd_i\in \cS(\cc_i)$ such that
  $\nor(\cd_i)\geq \nor(\cc_i)-r\cdot (k+1)$.
\end{Cor}

In other words:  If we assume that $\cK_i(n)$ is 
$(\card{I},r)$-decisive for all $i\in I$, then every
$\cc\in\cK_I(n)$ with $\nor(\cc)>1+r\cdot |I|$
is $(2,r\cdot(|I|+1))$-big.

In particular, we get pure decision for the finite product:
\begin{Cor}
  Set $\varphi(\mathord< n)\DEFEQ \prod_{i\in I}\prod_{m <\bfF_i (n)}\bfH_i(m)$,
  and $r(n)\DEFEQ 1/(n\varphi(\mathord< n))$.
  Assume that for all $i\in I$ and $n\in\omega$,
  $\cK_i(n)$ is
  $(\card{I},r(n))$-decisive
  and $r(n)$-halving.
  Then $\prod_{i\in I} \Qsu(\cK_i,\cS_i)$
  is $\omega\ho$-bounding and proper and preserves the size of the continuum.
  Under CH, $\prod_{i\in I} \Qsu(\cK_i,\cS_i)$ is $\al2$-cc and
  preserves all cardinals.
\end{Cor}

\begin{proof}
  $\prod_{i\in I} \Qsu(\cK_i,\cS_i)=\Qsu(\cK_I,\cS_I)$.
  $\cK_I(n)$  is $r(n)$-halving and $(2,r(n)\cdot (|I|+1))$-big according to
  Corollary~\ref{cor:multidimbignesssimple}.
  (Actually we get bigness only for creatures
  with norm bigger than $1+r\cdot |I|$ instead of $1$.)
  Now use Theorem~\ref{thm:onedim} and the Note following it.
  Note that $\prod_{i\in I} \Qsu(\cK_i,\cS_i)$ has size $2^\al0$.
\end{proof}

Remark:
    Decisiveness is quite costly: To be able to apply the
    last corollary, we will have to make the $n$-th level 
    much larger than levels before,
    i.e.
    \[
      \prod_{F_i(n)\leq m<F_i(n+1)}\bfH_i(m)\gg \prod_{j\in I}\prod_{m<F_i(n)}\bfH_j(m)
    \]
    for all $i\in I$.
    In our application this will have the effect that
    we can separate $(f,g)$ and $(f',g')$ only if their growth rates
    are considerably different. It is very likely that with a more careful
    and technically more complicated analysis one can construct
    forcings that can separate cardinal invariants for pairs that
    are not so far apart, but this would need other concepts than
    decisiveness.

\section{A variant of the countable support product}\label{sec:product}

We now define $P$, a variant of the countable support product of lim-inf
creature forcings. We want to end up with a forcing notion that also satisfies
fusion, pure decision and $\al2$-cc (under CH).  This will give preservation of
all cardinals. We will also need rapid reading of names.

Let $I$ be the index set of the product. We will use $\alpha$ and $\beta$ for
elements of $I$.
\begin{Asm}\label{asm:I}
  Fix a set $I$ and for every $\alpha\in I$,
  a creating pair $(\cK_\alpha,\cS_\alpha)$.
  We assume that for each $n$ there is an upper bound $m(n)$
  for $\card{\prod_{\bfF _\alpha(n)\leq i<\bfF _\alpha(n+1)}\bfH_\alpha(i)}$,
  and set \mbox{$\varphi(\mathord=n)\DEFEQ m(n)^n$},
  $\varphi(\mathord\leq n)\DEFEQ \prod_{m\leq n}\varphi(\mathord=m)$
  and $\varphi(\mathord< n)\DEFEQ \prod_{m<n}\varphi(\mathord=m)$.
\end{Asm}

We define the set $P$ in the following way:
\begin{Def}
  A condition $p$ in $P$ consists of a countable 
  subset $\dom(p)$ of $I$, of objects $p(\alpha,n)$ for $\alpha\in \dom(p)$, 
  $n\in\omega$, and of a function $\trunklg(p):\dom(p)\to \omega$  
  satisfying the following ($\alpha\in\dom(p)$): 
  \begin{itemize}
    \item
      If $n<\trunklg(p,\alpha)$, then
      $p(\alpha,n)\in \prod_{\bfF_\alpha(n)\leq i<\bfF_\alpha(n+1)}(\bfH_\alpha(i))$.
      \\
      $\bigcup_{n<\trunklg(p)}p(\alpha,n)$ is called trunk of $p$ at $\alpha$.
    \item
      If $n\geq \trunklg(p,\alpha)$, then $p(\alpha,n)\in \cK_\alpha(n)$
      and $\nor(p(\alpha,n))>0$.
    \item
      $\card{\supp(p,n)}<n$ for all $n>0$,
      where we set
      \[
        \supp(p,n)\DEFEQ \{\alpha\in\dom(p):\, \trunklg(p,\alpha)\leq n\}.
      \]
    \item
      Moreover, $\lim_{n\rightarrow\infty}(\card{\supp(p,n)}/n)=0$.
    \item
      $\lim_{n\rightarrow\infty}(\min(\{\nor(p(\alpha,n)):\, \alpha\in\supp(p,n)\}))= \infty$.
  \end{itemize}
\end{Def}
So in particular, 
for $\alpha\in\dom(p)$ the sequence $(p(\alpha,n))_{n\in\omega}$
is in $\Qsu(\cK_\alpha,\cS_\alpha)$.

Note that now there is an essential difference between a part $t$ of the
trunk and creature $\cc$ with $\val(\cc)=\{t\}$: The trunks do not prevent the
minimum of the norms at height $h$ to be large.

\begin{uRems}
  \begin{itemize}
    \item For the proof of Theorem~\ref{thm:ctbl}, we will additionally fix a
      function $\trunklg^\text{min}:I\to\omega$ and add the following
      requirement to the definition of $P$:
      \[
        \trunklg(p,\alpha)\geq \trunklg^\text{min}(\alpha).
      \]
      This does not change any of the following properties of $P$ (or their
      proofs).
    \item
      For the proof of Theorem~\ref{thm:uncountable}, we will define the forcing
      $R$ so that a condition $p$ picks for each $\alpha\in\dom(p)$ one of
      several possibilities for a creating pair $(\cK_\alpha,\cS_\alpha)$.  It
      turns out that this does not change anything either, apart from the fact
      that $R_\epsilon$ is not a complete subforcing of $R$ any
      more, i.e.\ Lemma~\ref{lem:Pjot} fails. Lemma~\ref{lem:cc}
      still holds but needs a new proof. The rest of the proofs still
      work without changes.
  \end{itemize}
\end{uRems}
 
As outlined, we have to modify the order usually used in the product:
\begin{Def}
  $q\leq p$ if
  \begin{itemize}
    \item $\dom(q)\supseteq \dom(p)$,
    \item if $\alpha\in\dom(p)$ and $n\in\omega$,
      then $q(\alpha,n)\in\cS(p(\alpha,n))$,
    \item $\trunklg(q,\alpha)= \trunklg(p,\alpha)$ for
          all but finitely many $\alpha\in\dom(p)$.
  \end{itemize}
\end{Def}
Note that $q\leq p$ implies that then $\trunklg(q,\alpha)\geq
\trunklg(p,\alpha)$ for all $\alpha\in\dom(p)$.

\begin{figure}[tb]
  \begin{center}
    \scalebox{0.25}{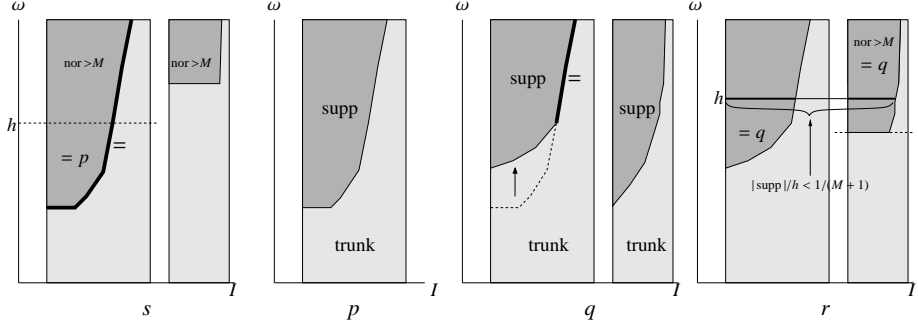}
    \caption{\label{fig:prod_pgtq}  $q\leq p$, $s\leq_M p$, $r\leq^\text{new}_M p$.}
  \end{center}
\end{figure}

Figure~\ref{fig:prod_pgtq}
shows one way to visualize $q\leq p$.

If $I$ is finite then
$P$ is just the product $\prod_{\alpha\in I}\Qsu(\cK_\alpha,\cS_\alpha)$.

For every $\alpha\in I$, $P$ adds a generic real $\n\eta_\alpha$, defined as
the union of the trunks of $p$ at $\alpha$ for $p$ in the generic filter. 
It is easy to see that $\n\eta_\alpha$ is forced to be different from
$\n\eta_\beta$ for $\alpha\neq \beta$. Once again, the sequence
$(\n\eta_\alpha)_{\alpha\in I}$ does not determine the generic filter.

Conditions with disjoint domains are compatible:
\begin{Lem}\label{lem:cc}
  (CH) $P$ is $\al2$-cc.
\end{Lem}
\begin{proof}
  Assume towards a contradiction that $A$ is an antichain of size
  $\al2$. Without loss of generality, $(\dom(a))_{a\in A}$
  forms a $\Delta$-system with root $u$.
  There are at most $2^{\al0}$ many possibilities for
  $a\restriction u$, so without loss of generality,
  $p\restriction u=q\restriction u$ for all $p,q\in A$.
  Then $p$ and $q$ are compatible:
  The function $x(n)=|\supp(p,n)\cup \supp(q,n)|/n$ converges
  to $0$. So there is an $h$ such that $x(m)<1$ for all $m\geq h$.
  Construct $r$ from $p\cup q$ by enlarging the (finitely many)
  trunks at $\supp(q,h)\cup \supp(p,h)$ to height $h$. Then $r\in P$ and
  $r\leq p,q$.
\end{proof}

\begin{Lem}\label{lem:Pjot}
  If $J\subseteq I$, then
  $P_{J}=\{p\in P:\, \dom(p)\subseteq J\}$ is a
  complete subforcing of $P$.
\end{Lem}

\begin{proof}
  If $p\in P$, then $p\restriction J\in P_J$, and
  $q\leq_P p$ implies $q\restriction J\leq_{P_J} p\restriction J$.
  So if $p\incomp_{P_J} q$, then $p\incomp_P q$.
  Also, $p\restriction J$ is a reduction of $p$:
  If $q\leq_{P_J} p\restriction J$, 
  then we can again enlarge finitely many stems
  of $q\cup p\restriction (I\setminus J)$ to get a condition $r\in P$ which
  is stronger than both $p$ and $q$.
\end{proof}

\begin{Def}
  \begin{itemize}
    \item $\prodval(p,\mathord< n)\DEFEQ \prod_{\alpha\in\dom(p)}\prod_{
      m<n}\val(p(\alpha,m))$. 
      The size of this set is at most $\varphi(\mathord<n)$.
      $\prodval(p,\mathord\leq n)\DEFEQ\prodval(p,<(n+1))$.
    \item If $w\subseteq \dom(p)$ and
      $t\in \prod_{\alpha\in w}\prod_{0\leq
          m<\bfF_\alpha(n)}\bfH_\alpha(m)$, then
          $p\wedge t$ is defined by
          \[
            (p\wedge t)(\alpha,m)=
        \begin{cases}
          t_\alpha\restriction [\bfF_\alpha(m),\bfF_\alpha(m+1)-1]  
             & \text{if }m<n\text{ and }\alpha\in w,\\
          p(\alpha,m) & \text{otherwise.}
        \end{cases}
      \]
      So $p\wedge t\in P$, and
          if $t\in \prodval(p,\mathord<n)$, then $p\wedge t\leq p$.
    \item If $\n\tau$ is a name of an ordinal, then
          $p$ $\mathord<n$-decides $\n\tau$, if  $p\wedge t$
          decides $\n\tau$ for all
          $t\in \prodval(p,\mathord< n)$. $p$ essentially
          decides $\n\tau$, if $p$ $\mathord<n$-decides $\n\tau$
          for some $n$.
  \end{itemize}
\end{Def}

As in the one-dimensional case we get:
\begin{Facts}\label{fact:variousfacts}
  \begin{enumerate}
    \item
      $\{p\wedge t:\, t\in\prodval(p,\mathord<n)\}$ is predense under $p$
      (for $p\in P$ and $n\in\omega$).
    \item $p\wedge t\, \forc\, \varphi$\text{ iff } $p\,\forc\,
        [(\forall \alpha\in\dom(t))\, t(\alpha)<\n\eta_\alpha\, \rightarrow\,
        \varphi]$.
    \item Assume that $q'$ is the result of replacing finitely
        many creatures $\cc$ of $q$ by creatures $\cd$ with
        $\val(\cd)\subseteq \val(\cc)$. Then $q'\leq^* q$.\footnote{%
        In other words: Assume that 
        $q,q'\in P$, $h\in\omega$,
        $\dom(q')=\dom(q)$, $q(\alpha,m)=q'(\alpha,m)$ for all
        $m\geq h$ and $\alpha\in\dom(q)$,  and
        $\val(q'_\alpha(m))\subseteq \val(q_\alpha(m))$ for all
        $m<h$ and $\alpha\in\dom(q)$.  Then $q'\leq^* q$.}
    \item If $q\leq p$ and $t\in\prodval(q,\mathord<n)$, then
        $t$ restricted to the domain of $p$ is in
        $\prodval(p,\mathord<n)$.\footnote{The same
          holds for $q\leq^* p$, apart from the fact that $\dom(p)$
          might not be a subset of $\dom(q)$. (Outside of $\dom(q)$,
          $p$ could consists of ``maximal creatures with no information''.)}
    \item If $q\leq p$, $t\in\prodval(q,\mathord<n)$, 
        and $s$ is the corresponding element in $\prodval(p,\mathord<n)$,
        then $q\wedge t\leq s\wedge p$.
    \item If $q'\leq q$ and 
      $q$ essentially decides $\n\tau$,
      then $q'$ essentially decides $\n\tau$.
    \item If $q\wedge t$ essentially decides $\n\tau$ for each
      $t\in\prodval(q,\mathord < n)$, then $q$ essentially decides $\n\tau$.
  \end{enumerate}
\end{Facts}

Recall that $\varphi(\mathord< n)$ 
is an upper bound for the number of possible sequences of trunks 
of height $n$ (cf.~\ref{asm:I}). 
\begin{Thm}\label{thm:prod}
  If $\cK_\alpha(n)$ is $(n,r(n))$-decisive
  and $r(n)$-halving for 
  $r(n)=1/(n^2\varphi(\mathord<n))$ and
  every $\alpha\in I$, $n\in\omega$,
  then $P$ is proper and $\omega\ho$-bounding.
  Assume $|I|\geq 2$ and set $\lambda=|I|^\al0$.
  Then $P$ forces $|I|\leq 2^\al0\leq \lambda$.
\end{Thm}

\begin{proof}
  The proof closely follows the one-dimensional case.
  We again prove pure decision and fusion, and the rest
  follows as in the proof of Theorem~\ref{thm:onedim}.
  (Note that $|P|=|I|^\al0$, and that $\n\eta_\alpha$ and $\n\eta_\beta$
  are forced to be different for $\alpha\neq\beta$.)

  So we have to define $\leq_M$: First we set
  $r\leq^\text{new}_{M} p$, if $r\leq p$, and
  \begin{itemize}
    \item if $n\in\omega$ and $\alpha\in\supp(r,n)\setminus \dom(p)$,
      then $n>M$, $\card{\supp(r,n)}/n\leq  1/(M+1)$, and
      $\nor(r(\alpha,n))> M$.
  \end{itemize}

  Assume that $M\in\omega$ and $q\leq p$. By  extending finitely many trunks in
  $q$ at positions $\alpha\notin\dom(p)$, we get an $r\leq q$ such that
  \begin{equation}\label{eq:trex}
    r\leq^\text{new}_{M} p\text{ and }
    r(\alpha,n)=q(\alpha,n)\text{ for }\alpha\in\dom(p)
  \end{equation}
  (cf. Figure~\ref{fig:prod_pgtq}).

  $s\leq^\text{old}_{M} p$, if $s\leq p$ and
  there is an $h\geq M$ such that
  for all $\alpha\in\dom(p)$,
  \begin{itemize}
    \item $\trunklg(s,\alpha)=\trunklg(p,\alpha)$,
    \item if $n<h$, then $s(\alpha,n)=p(\alpha,n)$,
    \item if $\alpha\in\supp(p,n)$ and $n\geq h$,
          then $\nor(s(\alpha,n))\geq M$.
  \end{itemize}

  $r\leq_{M} p$, if $r\leq^\text{new}_{M} p$ and $r\leq^\text{old}_{M} p$.

  By~\eqref{eq:trex} we get:
  \begin{equation}\label{eq:guto}
    \text{If $q\leq^\text{old}_{M} p$, then there is an
    $r\leq q$ such that $r\leq_{M} p$.}
  \end{equation}
  
  {\bf $\leq_n$ satisfies fusion:}\\
  Assume that $(p^m)_{m\in\omega}$ satisfies $p^{m+1}\leq_{m+1} p^{m}$.
  Define $q$ by $\dom(   q    )=\bigcup_{n\in\omega}\dom(p^n)$
  and $q_\alpha(n)=p^M_\alpha(n)$,
  where $M\geq n$ is minimal (or: arbitrary) such that $\alpha\in\dom(p^M)$.
  Then $q\in P$:
  Fix some $k$.
  Since $p^k\in P$, there is an $l$ such that 
  \begin{equation}\label{eq:xx7}
    \nor(p^k(\alpha,n))>k\text{ and }\card{\supp(p^k,n)}/n<1/(k+1)
    \text{ for all }n>l\text{ and }\alpha\in\supp(p^k,n).
  \end{equation}
  Since 
  \mbox{$p^{k+1}\leq_{k+1} p^k$, \eqref{eq:xx7}} holds 
  for $p^{k+1}$ as well, and for all $p^m$ with $m>k$, and
  therefore for $q$.

  So we just have to show pure decision:
  Fix $\n\tau$, a name of an ordinal.

  {\bf The basic construction $S(p,M)$:}
  \\
  Let $n$ be the minimal trunk-length of $p$, 
  i.e.\ $n=\min(\{\trunklg(p,\alpha):\, \alpha\in \dom(p)\})$.
  We will now define $S(p,M)\leq p$ for $M\in\omega$.

  Enumerate $\prodval(p,\mathord \leq n)$ as $s^0,\dots,s^{l-1}$.
  So $l\leq \varphi(\mathord=n)$.
  Set $p^{-1}\DEFEQ p$.
  Given $p^{k-1}$, define $p^k\leq p^{k-1}\wedge s^k$ and
  $h^{k}$ such that
  for all $\alpha\in \dom(p)$\footnote{we do not require
    anything for $\alpha\in \dom(p^k)\setminus\dom(p)$}
  \begin{itemize}
    \item
      $\trunklg(p^k,\alpha)=\trunklg(p^{k-1}\wedge s^k,\alpha)
      =\max(n+1,\trunklg(p,\alpha))$,
    \item
      if $n<m<h^k$,
      then $\nor(p^k(\alpha,m))\geq\nor(p^{k-1}(\alpha,m))-r(m)$,
    \item
      if $m\geq h^k$, then $\nor(p^k(\alpha,m))\geq M$,
  \end{itemize}
  and such that additionally one of the following two cases holds:
  \begin{description}
    \item[dec] $p^k$ essentially decides $\n\tau$, or
    \item[half] it is not possible to satisfy ``dec'' (for any choice
     of $h^k$),
     and $\dom(p^k)=\dom(p^{k-1})$ and
     $p^k(\alpha,m)=\chalf(p^{k-1}(\alpha,m))$ for all $m>n$ and
     $\alpha\in\supp(p^{k-1},m)$.
  \end{description}
  So we first try to find a $p^k$ satisfying ``dec'' (possibly
  with larger domain); if we fail we just halve each
  $p^{k-1}(\alpha,m)$.

  We construct $p^k$ for each $0\leq k<l$.
  This gives a function
  \[
    F:\prod_{\alpha\in\supp(p,n)}\val(p(\alpha,n))\fnto
     \{\text{dec},\text{half}\}.
  \]
  Each $\cK_\alpha(n)$ is $(n,r(n))$-decisive,
  and $\card{\supp(p,n)}< n$. So
  according to Corollary~\ref{cor:multidimbignesssimple} (for $k=n-1$)
  there are
  $\cd_\alpha\in\cS(p(\alpha,n))$ (for $\alpha\in\supp(p,n)$)
  such that $F\restriction \prod_{\alpha\in\supp(p,n)}\val(\cd_\alpha)$
  is constant and $\nor(\cd_\alpha)\geq \nor(p(\alpha,n))-n\cdot r(n)$.

  For this construction to work, we have to assume that the
  norms of all the creatures involved are big enough (so that
  we can apply bigness and halving). If
  this is not the case, $S(p,M)$ is undefined. Otherwise,
  we set $\dom(S(p,M))=\dom(p^{l-1})$
  and for $\alpha\in\dom(S(p,M))$
  \[
    S(p,M)(\alpha,m)=
    \begin{cases}
      p(\alpha,m) & \text{ if } m<n\text{ and }\alpha\in\dom(p),\\
      \cd_\alpha & \text{ if } m=n\text{ and }\alpha\in\dom(p),\\
      p^{l-1}(\alpha,m) & \text{ otherwise.}
    \end{cases}
  \]
  We call $q=S(p,M)$ halving, if the constant value of
  $F$ is ``half''.

  If $q$ is not halving, then $q$ essentially decides $\n\tau$:
  If $t\in\prodval(q,\mathord\leq n)$, then $t$ restricted to
  $\dom(p)$ is in
  $\prodval(p,\mathord\leq n)$, i.e.\ it is some
  $s^k$. 
  Then $q\wedge t\leq q\wedge s^k$, and
  $q\wedge s^k$
  is stronger than $p^k$, which essentially decides $\n\tau$. 
  Now use Facts~\ref{fact:variousfacts}(6,7).

  {\bf Some properties of $S(p,M)$:}
  \\
  If $q=S(p,M)$ is defined and $n$ the minimal trunk-length of
  $p$, then:
  \begin{gather}
    \label{eq:qamd}
    \parbox{0.8\columnwidth}{
      \raggedright
        $\nor(q(\alpha,n))\geq
        \nor(p(\alpha,n)) - n\cdot r(n)
        \text{ for } \alpha\in\supp(p,n).$
      }
    \\\label{eq:qbmd}
    \parbox{0.8\columnwidth}{
      \raggedright
        $\nor(q(\alpha,m))
        \geq\min(M,\nor(p(\alpha,m)))-\varphi(\mathord=n)\cdot r(m)
        $ for all $m>n$ and $\alpha\in\supp(p,m)$.
      }
    \\\label{eq:dreizweimd}
      \parbox{0.8\columnwidth}{
      \raggedright
        If $q$ is halving, then there is no $q'\leq q$
        essentially deciding $\n\tau$ such that
        \mbox{$\trunklg(q',\alpha)=\max(n+1,\trunklg(p,\alpha))$} for all
        $\alpha\in\dom(p)$.
      }
  \end{gather}
  To see~\eqref{eq:dreizweimd}, assume that $q'$ is a counterexample
  and that $h$ is such that
  $\nor(q'(\alpha,m))>M$ for all $m>h$ and $\alpha\in\supp(q',m)$.
  Let $t$ be in $\prodval(q',\mathord\leq n)$. $t$ restricted 
  to $\dom(p)$ is $s^k$ for some $k<l$.
  We know that 
  $p^k$ was constructed by halving each creature 
  of $p^{k-1}\wedge s^k$
  and that $q'\leq p^k$.
  We now define $r$: Set $\dom(r)=\dom(q')$.  If
  $m\leq h$ and
  $\alpha\in\supp(p,m)$, 
  we un-halve $q'(\alpha,m)$ to some 
  $\delta(\alpha,m)$ and set $r(\alpha,m)=\delta(\alpha,m)$.
  Otherwise we set $r(\alpha,m)=q'(\alpha,m)$.
  According to~\ref{fact:variousfacts}(3,6) $r$ essentially decides $\n\tau$.
  So we should have chosen $r$ instead of $p^k$, a contradiction.

  {\bf $S(p,M)$ essentially decides:}
  \\
  Assume that $M>3$, and
  that $\nor(p(\alpha,m))>3$ for all $m\in\omega$ and
  $\alpha\in\supp(p,m)$. We now show that  $S(p,M)$ exists
  and is not halving.

  Assume towards a contradiction that $S(p,M)$ is halving.
  Let $n$ be again the minimal trunk-length of $p$.
  We set $q^{n-1}=p$.
  Assume that for $k\geq n-1$, $q^k$ is
  already defined. We set
  $M_k=M+k+1-n$. (So $M_{n-1}=M$.)
  We define $q^{k+1}$ the following way:
  List  $\prodval(q^k,\mathord\leq k)$ as $s^0,\dots,s^{l-1}$.
  So $l\leq \varphi(\leq k)$. Set $r^{-1}\DEFEQ q^k$. Given
  $r^{i-1}$, set $r^i=S(r^{i-1}\wedge s^i,M_k)$ (if defined).
  Define $q^{k+1}$ to be $q^k$ up to $k$ and $r^{l-1}$ otherwise,
  and 
  additionally increase the stems outside $\dom(q^k)$
  to satisfy $q^{k+1}\leq^\text{new}_{M_k} q^k$.
  More formally: We pick some $h>M_k$, $h>k$  such that
  that $\nor(r^{l-1}(\alpha,m))>M_k$
  and $|\supp(r^{l-1},m)|/m<1/M_k$
  for all $m>h$ and $\alpha\in\supp(r^{l-1},m)$.
  For $\alpha\in \dom(r^{l-1})\setminus \dom(q^k)$
  and $m\leq h$, we pick some $t(\alpha,m)\in\val(r^{l-1}(\alpha,m))$.
  The we define $q^{k+1}$ by $\supp(q^{k+1})=\supp(r^{l-1})$ and
  \[
    q^{k+1}(\alpha,m)=
    \begin{cases}
      q^{k}(\alpha,m)& \text{ if }m\leq k\text{ and }\alpha\in\dom(q^k),\\
      r^{l-1}(\alpha,m)&  \text{ if }m>h\text{ or if }m>k\text{ and }\alpha\in\dom(q^k),\\
      t(\alpha,m)&  \text{ if }m\leq h\text{ and }\alpha\notin\dom(q^k).\\
    \end{cases}
  \]
  Note that $q^n$ is just $S(p,M)$ with some increased trunks 
  outside of $\dom(p)$.

  $q^{k+1}$ satisfies for $\alpha\in\dom(q^{k})$, $\beta\in\dom(q^{k+1})$:
  \begin{itemize}
    \item $q^{k+1}(\alpha,m)=q^k(\alpha,m)$ for $m\leq k$.
    \item $\nor(q^{k+1}(\alpha,k+1))\geq \nor(q^k(\alpha,k+1))-
          \varphi(\mathord\leq k)\cdot (k+1)\cdot r(k+1)$.
    \item $\nor(q^{k+1}(\alpha,m  ))\geq \min(M^k,\nor(q^k(\alpha,m  )))-
          \varphi(\mathord\leq k+1)\cdot r(m)$ for $m>k+1$.
    \item $\nor(q^{k+1}(\beta,m  ))\geq M^k$ if $\beta\in\supp(q^{k+1},m)\setminus \dom(q^k)$.
  \end{itemize}
  Iterating this $l$ many times, we get:
  \begin{align}
    \label{eq:qfmd}
    \nor(q^{k+l}(\alpha,m  ))\geq & \min(M^k,\nor(q^k(\alpha,m  )))-
              \min(l,m-k)\cdot \varphi(\mathord<m)\cdot m \cdot r(m),
    \\\intertext{so according to the definition of $r(m)$ we get}
    \label{eq:qgmd}
    \nor(q^{k+l}(\alpha,m  ))\geq & \min(M^k,\nor(q^k(\alpha,m  )))-1.
  \end{align}
  This shows, as in the one-dimensional case,
  that each $q^m$ is defined, and that $q^\omega$ is a condition in $P$,
  where we define
  $q^\omega$ by $\dom(q^\omega)=\bigcup_{k\in\omega} q^k$,
  and $q^\omega(\alpha,m)=q^k(\alpha,m)$, where $k$ is the minimal (or: some)
  $k\geq m$ such that $\alpha\in\dom(q^k)$. Just as for~\eqref{eq:prop7},
  we can generalize~\eqref{eq:dreizweimd} by induction and get:
  \begin{equation}\label{eq:prop7md}
    \parbox{0.8\columnwidth}{%
      \raggedright
        There is no $q'\leq q^{k}$ essentially deciding $\n\tau$ such that
        \mbox{$\trunklg(q',\alpha)=\max(k+1,\trunklg(q^{k},\alpha))$} for all
        $\alpha\in\dom(q^{k})$.
     }
  \end{equation}
  But there
  is a $q'\leq q^\omega$ deciding $\n\tau$. This implies that the
  trunk-lengths of $q'$ and of $q^\omega$ are the same on almost
  all elements of the domain of $q^\omega$. So by increasing 
  finitely many trunks of $q'$, we can assume that 
  $\trunklg(q',\alpha)=\max(k+1,\trunklg(q^\omega,\alpha))$
  for some $k$. So $q'\leq q^{k}$ decides 
  $\n\tau$, a contradiction to~\eqref{eq:prop7md}.\footnote{So
    this step in the proof is the reason that we had to redefine
    $\leq$.}

  {\bf Pure decision:}
  \\
  Given $p$ and $M$, we find an $h>M+6$ such that
  $\nor(p(\alpha,m))>M+6$ for all $m\geq h$ and $\alpha\in\supp(p,m)$.
  Enumerate $\prodval(p,\mathord \leq h-1)$ as
  $\{s^1,\dots, s^l\}$. As above,
  set $p^0=p$, $p^{k+1}=S(p^k\wedge s^k,M+6)$, and define
  $q$ by $q(\alpha,m)=p(\alpha,m)$ for $m<h$ and $\alpha\in\dom(p)$,
  and by $q(\alpha,m)=p^{l-1}(\alpha,m)$ otherwise.
  Then $q\leq^{\text{old}}_M p$ essentially decides $\n \tau$,
  and according to~\eqref{eq:guto} we
  find a $q'\leq q$ such that $q\leq_M p$.
\end{proof}

As already mentioned, only the growth rate of $r(n)$ is relevant.
Since we are dealing with decisive creatures, we can 
increase bigness even exponentially (in $n$) while decreasing the norms 
by a constant factor (cf.\ Corollary~\ref{cor:multidimbignesssimple}).
We use this for the following version of rapid reading.
Again, we set $\EXP(B,n,0)=B$ and $\EXP(B,n,k+1)=2^{\EXP(B,n,k)^n}$;
and we define $r$, $\varphi$ as in the previous theorem.

\begin{Thm}\label{thm:rapid} Assume that
  \begin{itemize} 
    \item $\delta\in\omega$,
    \item $g:\omega\fnto\omega$ is monotonously increasing,
    \item $\cK_\alpha(n)$ is  $(g(n),r(n))$-big,
      $(n,r(n))$-decisive and $r(n)$-halving
      for all $\alpha\in I$, $n\in\omega$,
    \item 
      $\n \nu(n)$ is a $P$-name and $p\in P$ forces that
      $\n \nu(n)<\EXP(g(n),n,n\cdot \delta)$ for all $n$.
    \end{itemize}
  Then there is a $q\leq p$ which 
  $\mathord<n$-decides $\n \nu(n)$ for all $n$ .
\end{Thm}

\begin{proof}
  We make the same modification to the previous proof as 
  in the one-dimensional case:

  {\em The basic construction $S(p,l,M)$:}
  We again assume that $n$ is the minimal length 
  of the trunks in $p$, and use the notation
  $S(p,l,M)$ (for $l\leq n$)
  for the same construction as $S(p,M)$, where 
  we set $\n\tau=\n\nu(l)$, and
  instead of trying to {\em essentially} decide $\n\tau$,
  we try to {\em decide} it.

  So 
  instead of the two cases ``dec'' and ``half'', we get 
  $\EXP(g(n),n,n\cdot \delta)+1$ many cases: one for each potential value
  of $\n\nu(n)$, and
  (if none of these cases can be satisfied) ``half''.
  So the number of possible cases is less than $\EXP(g(n),n,n\cdot (\delta+1))$.
  We use Corollary~\ref{cor:biggerbigness} to find
  successors $q(\alpha,m)$ of $p(\alpha,m)$
  with $F$-homogeneous product.
  This decreases the norm by at most $r(n)\cdot (n(\delta+1)+n)$,
  i.e.\ by $n\cdot (\delta+2)\cdot r(n)$.

  {\em Some properties of  $S(p,l,M)$:}
  So instead of~\eqref{eq:qamd} we get
  \[
        \nor(q(\alpha,n))\geq
        \nor(p(\alpha,n)) - n(\delta+2)\cdot r(n)
        \text{ for } \alpha\in\supp(p,n).
  \]
  There is no change to~\eqref{eq:qbmd}, and in~\eqref{eq:dreizweimd}
  we replace ``essentially deciding $\n\tau$'' with
  ``deciding $\n\nu(l)$''.

  {\em $S(p,l,M)$ decides:} We again construct $q^k$, each time trying to
  decide $\n\tau=g(l)$ (independently of $k$). 
  Instead of~\eqref{eq:qfmd}, we now get:
  \[
    \nor(q^{k+l}(\alpha,m  ))\geq \min(M^k,\nor(q^k(\alpha,m  )))-
    \min(l,m-k)\cdot \varphi(\mathord<m)\cdot m(\delta+2) \cdot r(m),
  \]
  and $r(m)=1/(m^2\varphi(\mathord<m))$.
  So 
  \[
    \min(l,m-k)\cdot \varphi(\mathord<m)\cdot m\cdot(\delta+2) \cdot r(m)
    \leq  m^2 \cdot \varphi(\mathord<m)\cdot r(m)\cdot (\delta+2)\leq \delta+2.
  \]
  So if we assume that
  \begin{equation}\label{eq:tmpgurke}
    \nor(p(\alpha,m))>\delta+2\text{ for all }m\in\omega\text{ and }
    \alpha\in\supp(p,m), 
 \end{equation}
  then again each $q^k$ (and $q^\omega$) is
  defined, and we get~\eqref{eq:prop7md} 
  for ``deciding $\n\nu(l)$''
  instead of ``essentially deciding $\n\tau$''.
  But there is some $q'\leq q^\omega$ deciding $\n\nu(l)$,
  a contradiction.

  So far we know the following:
  \begin{equation}\label{eq:gulag}
    \parbox{0,8\columnwidth}{
      \raggedright If $n$ is minimal trunk-length of $p$,
      if $p$ satisfies~\eqref{eq:tmpgurke}, and if $M>2(\delta+2)$, 
      then $S(p,n,M)$ exists and decides $\nu(n)$.}
  \end{equation}
  {\em Rapid reading:}
  Instead of the part on {\em pure decision}, we again proceed as follows:
  Fix $p\in P$ and $M>\delta+2$.  We can assume that 
  $p$ satisfies~\eqref{eq:tmpgurke},
  even for $2(\delta+2)$ instead of $\delta+2$ (just increase finitely many of
  the trunks).
  We set $k_0$ to be the minimal trunk-length of $p$, and
  $q^{k_0}=p$.
  We now construct $q^{k+1}$ and $q^\omega$ just as above, but this 
  time using 
  \begin{equation*}
    r^{i}=S(r^{i-1}\wedge s^i,k+1,M_k).
  \end{equation*} 
  I.e.\ we try to decide $\n\nu(k+1)$.
  Each $r^{i}(\alpha,n)$ has sufficient norm, and so 
  according to~\eqref{eq:gulag} $r^{i}$ (which has trunk-length $k+1$)
  decides $\n\nu(k+1)$. This implies that $q^{k+1}$ (and therefore
  $q^\omega$ as well) $\mathord\leq k$-decides $\n\nu(k+1)$.
\end{proof}

The rest of this section can safely be ignored: We describe how we
end up with our particular definition of the product. We want to find a
construction, similar to the countable support product, so that we can
generalize the pure decision proof of Section~\ref{sec:puredecondim}:
\begin{figure}[tb]
  \begin{center}
    \scalebox{0.25}{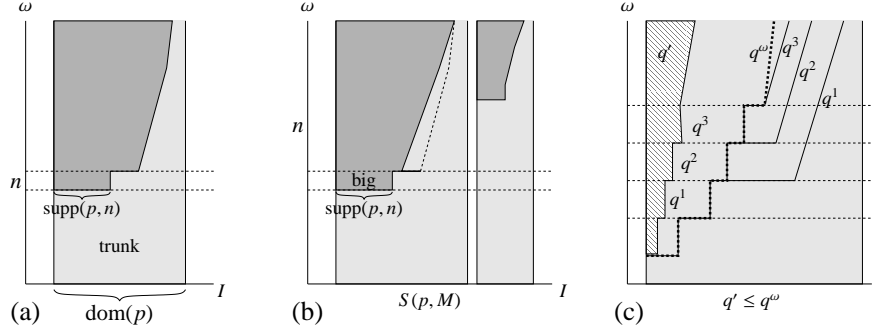}
    \caption{\label{fig:motivation}
      (a) A condition $p$ in $P$: $\dom(p)\subseteq I$ is 
      countable, at height $n$ there are less than $n$ many creatures.
      (b) The construction analog to $S(p,M)$.
      (c) We have to redefine $\leq$.
      }
  \end{center}
\end{figure}

\begin{itemize}
  \item To get $\al2$-cc, the support of the product can be at most
    countable.  For fusion, we have to allow at least countable support.
  \item A condition $p$ is a sequence
    $(p(\alpha,n))_{n\in\omega,\alpha\in\dom(p)}$. At each index $\alpha$,
    $p$ has a trunk, and above that $p(\alpha,n)$ is a creature 
    in $\cK_\alpha(n)$.
  \item To construct $S(p,M)$, we will set $n$ to be the minimal
    height of any stem of $p$. For each combination for
    values at height $n$ we get ``dec'' or ``half''.
    We want to use decisiveness to get homogeneous successors.
    For this we need that at height $n$,
    there are e.g.\ less than $n$ many
    creatures, and that $\cK(n)$
    is sufficiently decisive and big with
    respect to $n$. So we will generally assume that at
    each height $h$, there are less than $h$ many creatures,
    the rest is trunks, cf.\ Figure~\ref{fig:motivation}(a).
  \item 
    In the same construction step
    we also have to assume that 
    each of the creatures at height $n$ has sufficient norm.
    So we will not just require that for each $\alpha\in I$
    the norms of $p(\alpha,h)$ go to infinity, but that 
    the minimum of all the norms at height $h$ go to infinity.
  \item When we set $q=S(p,M)$ and are in the case ``half'', instead
    of~\eqref{eq:dreizwei}: 
    ``no $q'\leq q$ with trunk-length $n+1$
    essentially decides $\n\tau$'', we naturally get 
    ``no $q'\leq q$ essentially decides $\n\tau$, if
    the trunk-length at $\alpha$
    is the maximum of $n+1$ and the trunk-length of $p$ at $\alpha$.''
  \item We now assume towards a contradiction that $q=S(p,M)$ 
    is halving.
    We iterate the construction for all heights, get
    $q^\omega$, and find some  $q'\leq q^\omega$ essentially
    deciding $\n\tau$. However, this is not a contradiction:
    $q'$ could just have a longer trunk at each $\alpha$,
    cf.\ Figure~\ref{fig:motivation}(c).
  \item To fix this problem we redefine 
    $q\leq p$: We require that the trunk-lengths of $q$
    are (on the common domain) almost always equal to
    those of $p$, cf.\ Figure~\ref{fig:prod_pgtq}.
  \item Once we redefine $q\leq p$ this way, and additionally
    require that at level $h$ there are less than $h$ many
    creatures, we could end up with a
    condition whose domain cannot be enlarged any more 
    (since there already are maximally, i.e.\ $h-1$, many 
    creatures at each level $h$). We fix this by adding e.g.\ the requirement
    that the number of creatures at level
    $h$ divided by $h$ converges to $0$.
\end{itemize}

\section{A decisive creature with bigness and halving}\label{sec:example}

In this section, we construct decisive creatures with halving.

We use $\bfF (n)\DEFEQ n$ for all $n$, i.e.\ the $n$-creatures live on the
singleton $\{n\}$.

\begin{Lem}\label{lem:decexists}
  Assume that $n$ and $B$ are natural numbers, and that $0<r<1$.
  Then there is a natural number $\Psi(n,B,r)$ so that we can
  set $\bfH(n)=\Psi(n,B,r)$ and find 
  $r$-halving, $(B,r)$-big and $(n,r)$-decisive
  $n$-creatures 
  $(\cK(n),\cS)$ such that
  $\nor(\cc)>n$ for some $\cc\in\cK(n)$.
\end{Lem}

\begin{uRems}
  \begin{itemize}
    \item
      Without the last requirement the lemma is trivial, just assume that
      $\nor(\cc)<1$ for all $\cc\in\cK(n)$, and read the definitions of
      halving, big and decisive. 
    \item If such $(\cK(n),\cS)$ exists for some $\bfH(n)$, then it exists for
      every larger $\bfH(n)$ as well.
  \end{itemize}
\end{uRems}

The rest of this section consists of the proof of the lemma.
This proof is not needed in the rest of the paper.

We set $\flarge(m)=2^{2^{m^2}}$ and $a\DEFEQ 2^{\frac{1}{r}}$.  So
$\log_a(2)=r$.

\subsection*{The pre-norms} 
\begin{Lem}
  There is a $J\in\omega$ and a function $\ns$ on
  the powerset of $J$ such that the following holds:
  \begin{enumerate}
    \item $\ns$ is monotonous, i.e.\ %
      $u_1\subseteq u_2$ implies $\ns(u_1)\leq \ns(u_2)$.
    \item $\ns(\emptyset)=0$, and $\ns(J)\geq a^{n+1}$.
    \item If $\ns(u)=k+1$ then there is an $M\in\omega$ and a sequence
      $0=j_0<j_1<\dots<j_M$ such that $M\geq \max(B,\flarge(j_1+n))$ and
      $\ns(u\cap [j_i,j_{i+1}-1])\geq k$ for all $i\in M$.
  \end{enumerate}
\end{Lem}
\begin{proof}
  For finite subsets $u$ of $\omega$ define $\ns(u)\geq k$ by induction
  on $k$: For all $u$ set
  $\ns(u)\geq 0$, and $\ns(u)\geq 1$ iff $u$ is nonempty.
  For $k\geq 1$, we set $\ns(u)\geq k+1$ iff (3) as above holds. We show by
  induction on $k$ that for every $a\in\omega$ there is a $b\in\omega$ such that
  $\ns([a,b-1])=k$: Assume this is true for $k$. Given $a=j_0$,
  let $j_1$ be minimal such that $\ns([j_0,j_1-1])=k$.
  For every $i<\max(B,\flarge(j_1+n))$, find the minimal
  $j_{i+1}$ such that $\ns([j_i,j_{i+1}-1])=k$.
  Then $\ns([j_0,j_M-1])=k+1$.
  So we can pick $J$ such that $\ns([0,J-1])=a^{n+1}$.
\end{proof}
We set $\Psi(n,B,r)=\bfH(n)=2^J$. 
For a subset $c$ of $\bfH(n)$, we set
\[
  \prenor(c)\DEFEQ\max\{\ns(u):\, u\subseteq J, c\restriction u=2^u\},
\]
where $c\restriction u$ is $\{b\restriction u:\, b\in c\}$.
So $d\subseteq c$ implies $\prenor(d)\leq \prenor(c)$.

\begin{Lem}\label{lem:gurke2}
  Assume that $M\in\omega$, $J$ a set,
  $u\subseteq J$, $c\subseteq 2^{J}$, $c\restriction
  u=2^u$, $c=\bigcup_{i\in M}c_i$, and that $u_i$ ($i\in M$) are pairwise
  disjoint subsets of $u$.
  Then $2^{u_i}=c_i\restriction u_i$  for some $i\in M$.
\end{Lem}

\begin{proof}
  Otherwise, for all $i\in M$ there is an $a_i\in 2^{u_i}\setminus
  (c_i\restriction u_i)$.  Let $b\in 2^u$ contain the concatenation of these
  $a_i$. Then $b\in c\restriction u$, so $b\in c_i\restriction u$ for some $i\in
  M$, and $a_i\in c_i\restriction u_i$, a contradiction.
\end{proof}

\subsection*{The creatures}
An $n$-creature $\cc$ is a pair $(c,k)$ such that
$c\subseteq \bfH(n)$, $k\in\omega$ and $k\leq \prenor(c)-1$.
$\nor(\cc)$ is determined from $(c,k)$ by
\[
  \nor(c,k)\DEFEQ \log_{a}(\prenor(c)-k).
\]
For $n$-creatures $\cc\cong(c,k)$ and $\cd\cong(d,k')$ we define
\[
  (d,k') \in \cS(c,k)\text{ if }d\subseteq c\text{ and }k'\geq k.
\]

We now show that these creatures satisfy our requirements:
\begin{proof}[Proof of Lemma~\ref{lem:decexists}]
  It is clear that norms can be bigger than $n$:
  \[ 
    \nor(\bfH(n),0)=\log_a(\prenor(\bfH(n)))=\log_a(\ns(J))\geq \log_a(a^{n+1})=n+1.
  \]
  \subsubsection*{Halving}
  Assume $\nor(\cc)>1$, i.e.\ $\prenor(c)-k>a>2$.
  We define
  \[
    \chalf(c,k)\DEFEQ (c,k+\lfloor{(\prenor(c)-k)/2}\rfloor).
  \]
  Note that $\log_{a}(\lceil{(\prenor(c)-k)/2}\rceil)\geq
  \nor(c,k)-\log_{a}(2)= \nor(c,k)-r$.
  So
  \[
    \nor(\chalf(c,k))=
    \log_a(\prenor(c)-k-\lfloor{(\prenor(c)-k)/2}\rfloor)\geq
    \nor(c,k)-r.
  \]
  If $(d,k')\in \cS(\chalf(c,k))$ and $\nor(d,k')>0$, then
  \[
    \prenor(d)\geq
    k'+1\geq
    k+\lfloor{(\prenor(c)-k)/2}\rfloor+1,
  \]
  and we can un-halve $(d,k')$ to $(d,k)\in\cS(c,k)$:
  \[
    \nor(d,k)=
    \log_a(\prenor(d)-k)\geq
    \log_a(\lfloor{(\prenor(c)-k)/2}\rfloor+1)\geq
    \nor(c,k)-r,
  \]
  and $\val(d,k)=\val(d,k')=d$.

  \subsubsection*{Bigness}
  Let $(c,l)$ be an $n$-creature and $\nor(c,l)=x+r\geq r$.
  Let $u\subseteq J$ witness $\prenor(c)=a^{x+r}+l=2a^x+l$.
  So there is  an increasing sequence $(j_i)_{i\in M+1}$
  such that $c\restriction u=2^u$ and
  \begin{gather*}
    M\geq \max(B,\flarge(j_1+n)) \text{, and }
    \\
    \ns(u\cap [j_i,j_{i+1}-1])\geq 2a^x+l-1\geq a^x+l
    \text{ for all }i\in M. 
  \end{gather*}
  (If $x>0$, the last inequality is strict.)

  Take any $F:c\fnto M$. Then $c=\bigcup_{i\in M}F^{-1}\{i\}$.
  We set $u_i\DEFEQ u\cap[j_i,j_{i+1}-1]$ for $i\in M$. 
  According to Lemma~\ref{lem:gurke2} there is
  an $i\in M$ such that $F^{-1}\{i\}\restriction u_i=2^{u_i}$.
  We set $d\DEFEQ F^{-1}\{i\}\subseteq c$.
  Since $\ns(u_i)\geq a^x+l$ and $d\restriction u_i=2^{u_i}$,
  $\nor(d,l)\geq \log_a(a^x)= x=\nor(c,l)-r$.
  This shows that $(c,l)$ is $(M,r)$-big, and in particular $(B,r)$-big.

  \subsubsection*{Decisiveness}
  Pick $(c,l)\in \cK(n)$ such that  $\nor(c,l)=x+r\geq r$.
  As above there is a witness $u\subseteq J$, $M$
  and $(j_i)_{i\in M+1}$.
  Set $u^-\DEFEQ u\cap [j_0,j_{1}-1]$.
  Let $d^-\subseteq c$ contain
  for every $a\in 2^{u^-}$ exactly one $b\in c$
  such that $b\restriction u^-=a$. Then
  $\card{d}\leq 2^{j_1}\EQDEF K$ and (as above)
  $\nor(d^-,l)\geq \nor(c,l)-r$. So $(d^-,l)$ is a $K$-small
  successor of $(c,l)$.

  It remains to be shown that there is a $K$-big successor $(d^+,l)$.

  Let $F:c\fnto 2^{j_1}<M$ map $b$ to $b\restriction j_1$.
  So as above there is an $i<M$ such that
  $F^{-1}\{i\}\restriction u_i=2^{u_i}$ for $u_i\DEFEQ u\cap[j_i,j_{i+1}-1]$.
  Obviously $i\neq 0$.
  Set $d^+\DEFEQ F^{-1}\{i\}$.
  Pick any $(d',l')\in \cS(d^+,l)$ with norm bigger than 1.
  Let $\prenor(d')$ be witnessed by $u',M',(j'_i)_{i\leq M'}$.
  Then $u'\cap j_1=\emptyset$
  (since every $b\in d'$ has the same $b\restriction j_1$).
  So $j'_1>j_1$, and
  (by the same argument as above) $(d',l')$ is
  $(\flarge(j_1+n),r)$-big. This finishes the proof,
  since
  \[
    \flarge(j_1+n)=2^{2^{(j_1+n)^2}}\geq
    2^{2^{j_1\cdot n}}=2^{{(2^{j_1})}^n}=2^{K^n}.\qedhere
  \]
\end{proof}

\section{countably many cardinal invariants}\label{sec:countable}

Recall that $\myc_{f,g}$ and $\mycfa_{f,g}$ were defined in the introduction.

\newcommand{\fmax}{f_\text{max}}

In the previous section, we defined $\Psi(n,M,r)$
for $r>0$ and $n,M\in\omega$.
We can now specify the requirements we need for 
Theorem~\ref{thm:ctbl}:
\begin{Asm}\label{asm:sufficientlydifferent}
  $(f_\epsilon,g_\epsilon)_{\epsilon\in\omega}$ 
  is a sequence of functions from $\omega$ to $\omega$.
  $\fmax$ is such that $f_\epsilon(m)\leq \fmax(m)$ for 
  all $\epsilon\in\omega$.
  We set
  \[
    \varphi(\mathord=m)\DEFEQ\fmax(m)^m,
    \quad
    \varphi(\mathord<n)\DEFEQ\prod_{m<n}\varphi(\mathord=m)
    \quad
    r(n)\DEFEQ \frac{1}{n^2\varphi(\mathord<n)},
  \]
  and assume:
  \begin{itemize}
    \item If $\epsilon\neq \epsilon'$, then there is an
      $n$ such that $f_\epsilon(m)\neq f_{\epsilon'}(m)$
      for all $m>n$.
    \item $f_{\epsilon}(m)\gg g_{\epsilon}(m)$ for all $\epsilon,m$; 
      more precisely $f_{\epsilon}(m)\geq \Psi(m,g_\epsilon(m),r(m))$.
    \item If $f_\epsilon(m) > f_{\epsilon'}(m)$,
      then $g_{\epsilon}(m)\gg f_{\epsilon'}(m)$; more precisely
      $\varphi(\mathord<m)f_{\epsilon'}(m)^m<g_\epsilon(m)$. 
    \item $g_\epsilon(m)>\varphi(\mathord<m)$.
    \item $g_{\epsilon}(m+1)\geq \fmax(m)$
      for all $\epsilon,m\in\omega$.
  \end{itemize}
\end{Asm}

The assumption states more or less that the 
$f_\epsilon,g_\epsilon$ have sufficiently different 
growth rates, and that each level is much bigger than
the previous levels.
If is clear that we can construct such sequences (by induction).

\begin{Thm}\label{thm:main}
  Assume CH.
  Choose for all $\epsilon\in\omega$ a cardinal $\kappa_\epsilon$ such that
  $\kappa_\epsilon=\kappa_\epsilon^\al0$.
  Let $(f_\epsilon,g_\epsilon)_{\epsilon\in \omega}$ be as above.
  Then there is a proper, $\al2$-cc, $\omega^\omega$-bounding
  partial order $P$ which preserves cardinals and forces
  that $\myc_{f_\epsilon,g_\epsilon}=\mycfa_{f_\epsilon,g_\epsilon}=
  \kappa_\epsilon$ for all $\epsilon\in \omega$.
\end{Thm}

Let $I$ be the disjoint union of $I_\epsilon$ ($\epsilon\in\omega$) such that
each $I_\epsilon$ has size $\kappa_\epsilon$ and is disjoint to $\omega$.

We will use $\epsilon, \epsilon', \epsilon_1, \dots$ for the cardinal
invariants (i.e.\ for elements of $\omega$), and $\alpha, \beta, \dots$ for
elements of $I$. $I$ will be the  index set of the product.

So according to the definition of $\Psi$,
we can choose for each $\epsilon,n\in\omega$ a
creating pair $(\cK_\epsilon(n), \cS_\epsilon)$
satisfying the following:
\begin{itemize}
  \item $\bfF_\epsilon (n)=n$,
  \item $\bfH_\epsilon(n)=f_\epsilon(n)$,
  \item $\cK_\epsilon(n)$ is
    $(g_\epsilon(n),r(n))$-big, $r(n)$-halving and
    $(n,r(n))$-decisive.
\end{itemize}

For every $\alpha\in I_\epsilon$ and $n\in\omega$, we set $\cK_\alpha(n)\DEFEQ
\cK_{\epsilon}(n)$, $f_\alpha \DEFEQ f_{\epsilon}$ and $g_\alpha \DEFEQ
g_{\epsilon}$ and we set
$\trunklg^\text{min}(\alpha)$ to be the minimal $n$ such that
$f_{\epsilon'}(m)\neq f_\epsilon(m)$ for all $\epsilon'<\epsilon$.

$P$ is the forcing notion defined in Section~\ref{sec:product}, where we
additionally require 
\begin{itemize}
  \item $\trunklg(p,\alpha)\geq \trunklg^\text{min}(\alpha)$
    for all conditions $p$ and $\alpha\in\dom(p)$.
\end{itemize}
As already noted, this does
not change any of the results of Section~\ref{sec:product}.

Note that 
$\varphi(\mathord <n)$ and $r(n)$ are as in Theorem~\ref{thm:prod}, and
that we assume CH.
So we get:
\begin{Cor}\label{cor:variousbigness}
  \begin{enumerate}
    \item $P$ is proper and $\al2$-cc,
      $P$ has continuous reading of names, and preserves all cardinals.
    \item {\em (Separated support.)} 
       If $p\in P$, $\alpha,\beta\in\supp(p,n)$, $\alpha\in I_\epsilon$,
       $\beta\in I_{\epsilon'}$, and $\epsilon\neq\epsilon'$, then
       $f_\epsilon(n)\neq f_{\epsilon'}(n)$.
    \item {\em (Rapid reading.)} If $p\in P$ forces that $\n\eta$ is an
       $(f_\epsilon,g_\epsilon)$-slalom,
       or that $\n\eta(n)<f_\epsilon(n)$ for all $n$,
       then there is a $q\leq p$ which $\mathord\leq n$-decides
       $\n\eta(n)$ for all $n\in\omega$.
  \end{enumerate}
\end{Cor}

It also follows that $P_\epsilon\DEFEQ P_{I_\epsilon}$ is a complete
subforcing of $P$ and forces that the size of the continuum is 
$\kappa_\epsilon$.

\begin{proof}
  (1): Theorem~\ref{thm:prod} and Lemma~\ref{lem:cc}.
  (2): Assume that $\epsilon<\epsilon'$.
    $\trunklg(p,\beta)>\trunklg^\text{min}(\beta)$, i.e.\ $f_\epsilon(n)\neq
    f_{\epsilon'}(n)$. 
  (3) follows from~\ref{thm:rapid}: Set 
    $\delta=3$, $g(n)=\fmax(n-1)$ and
    $\n\nu(n)=\n\eta(n-1)$ for all $n$.
    Each $\cK_\epsilon(n)$ is  $(g_\epsilon(n),r(n))$-big
    for some $\epsilon$, $g_\epsilon(n)\geq \fmax(n-1)=g(n)$,
    and $p$ forces that there are at most
    $\fmax(n-1)^{\fmax(n-1)}<\EXP(g(n),n,3)$ many possible values 
    for $\n \nu(n)$.
    So there is a $q\leq p$ which $\mathord< n$-decides
    $\n\nu(n)=\n\eta(n-1)$ for all $n$.
\end{proof}

In the following two sections, we will show that $P$ forces
$\kappa_\epsilon\leq \myc_{f_\epsilon,g_\epsilon}$ and $
\mycfa_{f_\epsilon,g_\epsilon}\leq \kappa_\epsilon$.  This proves
Theorem~\ref{thm:ctbl}, since $\myc_{f,g}\leq  \mycfa_{f,g}$ for all $(f,g)$.

\section{$P_\epsilon$ adds a $\forall$-cover}\label{sec:forall}

\begin{Lem}\label{lem:forall}
  $P$ forces $\mycfa_{f_\epsilon,g_\epsilon}\leq \kappa_\epsilon$.
\end{Lem}

One nice way to formulate the proof is the following: $P_\epsilon$ is a complete
subforcing and forces $2^\al0=\kappa_\epsilon$. And in the $P$-extension
$V[G]$, the set of slaloms that are in the $P_\epsilon$-extension $V[G\cap
P_\epsilon]$ form a $(\forall,f_\epsilon,g_\epsilon)$-cover.

However, to be able to generalize the proof to the uncountable
case of Section~\ref{sec:uncountable}, we will not use the complete
subforcing. Instead we will use pure decision more explicitly.

\begin{proof}
  Let $p_0\in P$ and
  $\n r$ be a $P$-name for a real such that $\n r(n)<f_\epsilon(n)$
  for all $n$.
  We will show that
  \begin{equation}
     \label{eq:rapiddense}
     \parbox{0.8\textwidth}{
     There is a $q\leq p_0$ and a way
     to determine an $(f_\epsilon,g_\epsilon)$-slalom
     $\n S(n)$ from $\prodval(q,\mathord \leq n)$ restricted to $I_\epsilon$,
     such that $q$ forces $\n r(n)\in \n S(n)$ for all $n$.}
  \end{equation}

  More explicitly, we find a $q$ and a function $\text{eval}$ which
  assigns to each $t\restriction I_\epsilon$ for $t\in \prodval(q,\mathord \leq n)$
  a set $S^t(n)$ such that $S^t(n)\subseteq f_\epsilon(n)$,
  $|S^t(n)|\leq g_\epsilon(n)$ and such that $q$ forces the following:
  If $t$ is compatible with the generic filter, then $\n r(n)\in S^t(n)$.

  Assume that we can do this for all names $\n r$. 
  Note that there are only $\kappa_\epsilon$ many possible assignments
  as above: There are only $\kappa_\epsilon^\al0=\kappa_\epsilon$ many 
  possible sequences $q\restriction I_\epsilon$, and $2^\al0$ many 
  ways to continuously read  a real from $q\restriction I_\epsilon$.
  Each assignment, together with the $P$-generic filter, determines 
  a slalom $\n S$. Let $X$ be the set of all possible assignments.
  This corresponds to a $P$-name $Y$ of a family (of size $\kappa_\epsilon$)
  of $(f_\epsilon,g_\epsilon)$-slaloms, and according
  to~\eqref{eq:rapiddense}, the following holds in the $P$-extension:
  For every $\eta\in\prod_{n\in\omega} f_\epsilon(n)$ there is a slalom $\n S$
  in $Y$ covering $\eta$.
  This implies $\mycfa_{f_\epsilon,g_\epsilon}\leq
  \kappa_\epsilon$.
  
  So it remains to show~\eqref{eq:rapiddense}.
  First pick a $p\leq p_0$ rapidly reading $\n r$ as in~\ref{cor:variousbigness}(3), i.e.\ $p$ $\mathord\leq n$-decides
  $\n r(n)$ for all $n\in\omega$.
  We can assume that $\nor(p_\alpha(n))>3$ for all $\alpha\in\supp(p,n)$.
  We set $\dom(q)=\dom(p)$ and
  $\trunklg(q,\alpha)=\trunklg(p,\alpha)$, 
  and we will define $q(\alpha,m)$ (for all $\alpha\in\supp(p,m)$) as
  well as $\n S(m)$ by induction on $m$. We will 
  find $q(\alpha,m)\in \cS(p(\alpha,m))$ such that 
  the norm decreases by at most 2. Then $q$ 
  automatically is a valid condition in $P$ and stronger than $p$.

  Fix $m\in\omega$. Set $M\DEFEQ \supp(p,m)\cap I_\epsilon$.
  ($M$ stands for ``medium''.)
  According  to ``separated support''~\ref{cor:variousbigness}(2),
  \begin{equation}\label{eq:qwrr}
    \alpha\in \supp(p,m)\setminus I_\epsilon\text{ implies }
   f_{\alpha}(m)\neq f_\epsilon(m).
  \end{equation}
  So either $f_{\alpha}(m)<f_\epsilon(m)$, 
  in this case we set $\alpha\in S$ (for ``small'');
  or $f_{\alpha}(m)> f_\epsilon(m)$, then we set $\alpha\in L$
  (for ``large'').
  So $\supp(p,m)$ is partitioned into $S$, $M$ and $L$.
  We set $q(\alpha,m)=p(\alpha,m)$ for $\alpha\in S\cup M$.
 
  $p$ $\mathord\leq m$-decides $\n r(m)$, i.e.\ there is a function $F$
  that calculates $\n r(m)<f_\epsilon(m)$:
  \[
    F: \prodval(p,\mathord<m)
       \times
       \left(\prod_{\alpha\in S\cup M\cup L}\val(p_\alpha(m))\right)
       \fnto
       f_\epsilon (m).
  \]

  {\em Step 1:} 
  Assume $L$ is nonempty (otherwise continue with Step~2).
  \[
    \left|\prod_{\alpha\in S\cup M}\val(p_\alpha(m))\right|
    \leq \bfH_\epsilon(m)^{m-1}=f_\epsilon(m)^{m-1}.
  \]
  So we can rewrite $F$ as
  \[
    F': \prod_{\alpha\in L}\val(p_\alpha(m))\fnto
        f_\epsilon (m)^{\varphi(\mathord<m) f_\epsilon(m)^{m-1}}
        < f_\epsilon (m)^{ f_\epsilon (m)^m}.
  \]
  If we set $B=\min(\{g_\alpha(m):\, m\in L\})$, then
  $f_\epsilon(m)<B$, and $B^{B^m}<\EXP(B,m,3)$. 
  According to
  Corollary~\ref{cor:biggerbigness},
  there are $q(\alpha,m)\in \cS(p(\alpha,m)$ for $\alpha\in L$
  such that
  $F'$ restricted to $\prod_{\alpha\in L}\val(q(\alpha,m))$ is constant
  and $\nor(q(\alpha,m))>\nor(p(\alpha,m))-r(m)\cdot (m+3)$.
  This defines $q(\alpha,m)$ for $\alpha\in L$. So we now know
  $q(\alpha,m)$ for all $m$.

  {\em Step 2:} 
  So (modulo $q$)
  we have eliminated the dependence of $\n r(m)$ 
  on $L$, and are left with
  \[
    F: \prodval(q,\mathord<m)
       \times
       \left(\prod_{\alpha\in S\cup M}\val(q(\alpha,m))\right)
       \fnto
       f_\epsilon (m).
  \]
  We now define $\n S(m)$, more exactly the
  evaluation that maps $t\in\prodval(q,\mathord\leq m)\restriction I_\epsilon$
  to $S^t(m)$. So fix such a $t\in \prod_{\alpha\in M}\val(q(\alpha,m))$.

  $q\wedge t$ allows for at most
  $\varphi(\mathord<m)\cdot \prod_{\alpha\in S}\val(q(\alpha,m))$ 
  many possible values for $\n r(m)$.

  If $S$ is nonempty, let $\epsilon'$ be such that
  $f_{\epsilon'}(m)=\max\{f_\alpha(m):\, \alpha\in S\}$. Then 
  $\prod_{\alpha\in S}\val(p_\alpha(m))\leq f_{\epsilon'}(m)^m$.
  So we get 
  $\varphi(\mathord<m)\cdot f_{\epsilon'}(m)^m<g_\epsilon(m)$
  many possible values for $\n r(m)$.
  (If $S$ is empty,
  we just get $\varphi(\mathord<m)$ many possibilities.)
  So we can set $S^t(m)$ to be this set of possible values.
\end{proof}

\section{There is no small $\exists$-cover}\label{sec:exists}
\begin{Lem}\label{lem:exists} (CH) $P$ forces
  $\kappa_\epsilon \leq \myc_{f_\epsilon,g_\epsilon}$.
\end{Lem}

\begin{proof}
  Assume towards a contradiction that $p_0$ forces that $\n\SlFam$ is an
  $(\exists,f_\epsilon,g_\epsilon)$-cover, $\al1\leq \lambda<\kappa_\epsilon$
  and $\n\SlFam=\{\n S_i:\, i\in\lambda\}$.

  For every $i$, the set of $p'\leq p_0$ which 
  rapidly\footnote{as in Corollary~\ref{cor:variousbigness}(3)}
  reads $\n S_i$ is predense under $p_0$.
  Because of $\al2$-cc, we can find a set
  $D_i$ of such $p'$
  which is predense under $p_0$ and has size $\al1$.
  So
  \[
    J=\bigcup_{i\in\lambda,p'\in D_i}\dom(p')
  \]
  has size $\lambda$.
  Since $\card{I_\epsilon}=\kappa_\epsilon>\lambda$, there
  is a $\beta\in I_\epsilon\setminus J$. Fix this $\beta$.

  Let $p_1\leq p_0$ decide the $i$ such that  
  $\n \eta_\beta(n)\in \n S_i(n)$ for infinitely many $n$. We
  set $\n S\DEFEQ \n S_i$.
  We can assume $\beta\in\dom(p_1)$, so we have
  \begin{equation}\label{eq:bdomq}
    \beta\in\dom(p_1)\cap I_\epsilon\setminus J.
  \end{equation}
  Let $p\leq p_1$ be stronger than some $p'\in D_i$,
  and let $\nor(p(\alpha,m))>10$ for all $\alpha\in\supp(p,m)$.
  So modulo $p$, we can determine 
  the value of $\n S(n)$ from $t\restriction J$
  for $t\in \prodval(p,\mathord\leq n)$.\footnote{More formally: Let $X$ 
    be the set $\{t\restriction J:\, t\in\prodval(p,\mathord\leq n)\}$.
    For each $x\in X$ there is an $S^x_n$ such that $p$ forces:
    $(\forall \alpha\in J)\, x(\alpha)<\n\eta_\alpha\,\rightarrow\,
    \n S(n)=S^x_n$.}

  We will show towards a contradiction 
  that we can strengthen $p$ to a $q$ such that
  for all $n\geq \trunklg(p,\beta)$ the following holds:
  the generic $\n \eta_\beta(n)$ (which is  in
  $\val(q(\beta,n))$ and less than  $f_\epsilon(n)$)
  avoids every possible element of $\n S(n)$,
  (which is determined by $q(\alpha,m)$
  for $m\leq n$ and $\alpha\neq \beta$).
  In other words, we can make $\n\eta_\beta$ run away from
  $\n S$ at every height above the trunk.
  So $q$ forces that $\n \eta_\beta(n)\notin \n S(n)$
  for all $n\geq \trunklg(p,\beta)$, a contradiction.

  We set $\dom(q)=\dom(p)$, $\trunklg(q,\alpha)=\trunklg(p,\alpha)$,
  and define $q(\alpha,m)$ (for all $\alpha\in\supp(p,m)$)
  by induction on $m$. We will find a
  $q(\alpha,m)\in \cS(p(\alpha,m))$ so that the
  norm decreases by at most 2. This guarantees that
  $q$ is a condition in $P$ and stronger than $p$.

  Fix an $n\geq \trunklg(p,\beta)$.
  Set $A\DEFEQ\supp(p,n)$. So $\beta\in A$, and
  without loss of generality $\card{A}\geq 2$.
  According to the definition of $P$, $\card{A}<n$.

  Similarly to the previous section, we will partition 
  $A$ into the large indices $L$, the small ones $S$ and
  $\{\beta\}$. However, we cannot assume that 
  $A\cap I_\epsilon=\{\beta\}$, so the partition
  will not only be based on membership in $I_{\epsilon'}$,
  but has to be ``finer''. $\n S(n)$ only depends
  on $S\cup L$ (and the very small set $\prodval(p,\mathord<n)$).
  Again, we first use bigness to eliminate the dependence of $\n S(n)$
  on the large part. And the small part is sufficiently small
  so that $\n\eta_\beta(n)$ (i.e.\ $q(n,\beta)$) 
  avoids all the possible elements of $\n S(n)$.
  We now do this in more detail:

  \begin{figure}[tb]
    \newcommand{\mads}[1]{\ar[d]_-{#1\text{-small}}}
    \newcommand{\madb}[1]{\ar[d]^-{\mathord\geq#1\text{-big}}}
    \centerline{\xymatrix@R-1em{
      \cc^0_{\alpha_0}           \mads{K_0} & \cc^0_{\alpha_1}           \madb{K_0}& \dots &  \cc^0_{\alpha_m}=\cc^0_\beta\madb{K_0}    & \dots & \cc^0_{\alpha_{\card{A}-2}} \madb{K_0}        & \cc^0_{\alpha_{\card{A}-1}}\madb{K_0}\\
      \cd_{\alpha_0}                        & \cc^1_{\alpha_1}           \mads{K_1}     & \dots &  \cc^1_{\alpha_m}            \madb{K_1}    & \dots & \cc^1_{\alpha_{\card{A}-2}} \madb{K_1}        & \cc^1_{\alpha_{\card{A}-1}}\madb{K_1}\\
                                            & \cd_{\alpha_1}                            &       &  \vdots                      \madb{K_{m-1}}&       & \vdots                                             & \vdots                                    \\
                                            &                                           &       &  \cc^m_{\alpha_m}            \mads{K_m=K}       & \dots & \vdots                                             & \vdots                                    \\
                                            &                                           &       &  \cd_{\alpha_m}=\cd_{\beta}                     & \dots & \vdots                      \mads{K_{\card{A}-2}}  & \vdots                     \madb{K_{\card{A}-2}}    \\
                                            &                                           &       &                                                 &       & \cd_{\alpha_{\card{A}-2}}                                   & \cd_{\alpha_{\card{A}-1}}                 \\
    }}
    \caption{\label{fig:uiq}}
  \end{figure}
  Set $\cc^0_\alpha\DEFEQ p(\alpha,n)$ for $\alpha\in A$.
  Assume that for $l\geq 0$ we already have 
  a list $(\alpha_k)_{k<l}$ of elements of $A$
  and creatures
  $(\cc^l_\alpha)_{\alpha\in A\setminus\{\alpha_0,\dots,\alpha_{l-1}\}}$.
  Each $\cc^l_\alpha$ is $(K^l_\alpha,n,r(n))$-decisive
  for some $K^l_\alpha$.
  Set $K_l\DEFEQ\min(\{K^l_\alpha:\, \alpha\in A\setminus
  \{\alpha_0,\dots,\alpha_{l-1}\}\})$,
  and choose
  $\alpha_l$ such that $K_{\alpha_l}^l=K_l$.
  Let  $\cd_{\alpha_l}$ be  a $K_l$-small successor of $\cc^l_{\alpha_l}$.
  For $\alpha\in A\setminus \{\alpha_0,\dots,\alpha_{l}\}$, let
  $\cc^{l+1}_\alpha$ be 
  a $K_l$-big successor of $\cc^l_{\alpha}$.
  Cf.\ Figure~\ref{fig:uiq}.
  Iterate this construction $\card{A}-1$ times.
  So there remains one
  $\alpha$ that has not been listed as an $\alpha_l$, set 
  $\alpha_{\card{A}-1}=\alpha$ and
  $\cd_{\alpha_{\card{A}-1}}=\cc^{\card{A}-1}_\alpha$.
  
  Let $m$ be such that $\beta=\alpha_m$, and set 
  \[
    K\DEFEQ K_m,
    \quad
    S\DEFEQ \{\alpha_l:\, l<m\},
    \quad
    L\DEFEQ \{\alpha_l:\, l>m\}.
  \]
  So $A$ is partitioned into the three parts $\{\beta\}$, $S$ and
  $L$. We get:
  \begin{itemize}
    \item $\cd_\alpha\in \cS(p(\alpha,n))$,
      $\nor(\cd_\alpha)\geq \nor(p(\alpha,n))-(n-1)\cdot r(n)$.
    \item $\prod_{\alpha\in S}\card{\val(\cd_\alpha)}\leq K^{n-2}_{m-1}<K$.
    \item $\cd_\beta$ is hereditarily 
      $K_{m-1}$-big\footnote{even
       $2^{K^n_{m-1}}$-big. Provided of course that $S$ is nonempty,
       otherwise there is no $K_{m-1}$.}
       and $\card{\val(\cd_\beta)}\leq K$.
    \item If $\alpha\in L$, then $\cd_\alpha$ is hereditarily
      $K$-big.\footnote{even $2^{K^n}$-big.}

  \end{itemize}
  $J\cap \supp(p,n)\subseteq S\cup L$, so
  $\n S (n)$ is determined by
  $\prodval(p,\mathord<n)\times \prod_{\alpha\in S\cup
  L}\val(p(\alpha,n))$.
  We set $q(\alpha,m)=\cd_\alpha$ for all $\alpha\in S$.
  We also set $q(\beta,m)=\cd_m$ for now. (But we may further
  decrease $q(\beta,m)$ in Step~2.) We
  are only interested in the
  elements of $\n S (n)$ that are possible values of
  $\n\eta_\beta(n)$,
  in other words we are interested in
  $\n S (n)\cap \val(\cd_\beta)$.
  This part has size at most $K$.
  So we get a function
  \[
    F: \prodval(p,\mathord<n)\times
    \left(\prod_{\alpha\in S}\val(\cd_\alpha)\right)\times
    \left(\prod_{\alpha\in L}\val(\cd_\alpha)\right)\fnto
    \binom{K}{g_\epsilon(n)}.
  \]
  {\em Step 1:} Assume  $L$ is non-empty (otherwise
  continue with Step~2).
  Note that
  $\binom{K}{g_\epsilon(n)}\leq K^{g_\epsilon(n)}$ and
  $\varphi(\mathord<n)< g_\epsilon(n)<K$.
  So we can rewrite $F$ as
  \[
    F': \prod_{\alpha\in L}\val(\cd_\alpha)\fnto
     (K^K)^{K\times K}=K^{K^3}.
  \]
  Since $\cd_\alpha$ is decisive and (hereditarily)
  $K$-big for $\alpha\in L$
  and $\EXP(K,n,3)>K^{K^3}$,
  we can find 
  $F'$-homogeneous
  $q(\alpha,n)\in\cS(\cd_\alpha)$ for
  $\alpha\in L$ such that
  the norm decreases by at most $(n+1)\cdot r(n)$,
  cf.\ Corollary~\ref{cor:biggerbigness}.

  {\em Step 2:}
  So modulo $q$ we have eliminated $L$ and can rewrite $F$ as
  \[
    F: \prodval(p,\mathord<n)\times
    \left(\prod_{\alpha\in S}\val(\cd_\alpha)\right)\fnto
    \binom{K}{g_\epsilon(n)}.
  \]
  Let $X$ be the image of $F$ (i.e.\ the set of possible values
  of $\n S(n)\cap \val(\cd_\beta)$). $X$ has size at most
  $\varphi(\mathord<n)\cdot K^{n-2}_{m-1}<\EXP(K_{m-1},n,2)$.
  So according to~\ref{lem:increasebigness}(5), 
  we can strengthen $\cd_\beta$ to avoid $X$,
  decreasing the norm by at most $3\cdot r(n)$. 
\end{proof}

\section{Uncountably many invariants}\label{sec:uncountable}

We construct natural numbers
$(f_{n,l})_{n\in \omega,-1\leq l\leq n}$, and
$(g_{n,l})_{n\in \omega,0\leq l\leq n}$ 
so that $0=f_{0,-1}$ and (for $n,l\in\omega$) $f_{n+1,-1}=f_{n,n}$ and
$f_{n,l-1}< g_{n,l}< f_{n,l}$.
We set
$\fmax(n)=f_{n,n}$, 
$\varphi(\mathord =n)=\fmax(n)^n$,
$\varphi(\mathord <n)=\prod_{m<m}\varphi(\mathord =m)$
and $r(n)=1/(n^2\varphi(\mathord<n))$.
So  we get the following picture:\\
\setlength{\unitlength}{5mm}
\centerline{\begin{picture}(23,2)(-0.4,-1.2)
  \put( 0,0){\line(1,0){21}}
  \put( 0,-0.3){\line(0,1){0.6}}
  \put( 0,-1){\makebox(0,0){\footnotesize $0$}}
  \put( 3,-0.15){\line(0,1){0.3}}
  \put( 3,-0.7){\makebox(0,0){\footnotesize $g_{0,0}$}}
  \put( 6,-0.3){\line(0,1){0.6}}
  \put( 6,-1){\makebox(0,0){\footnotesize $f_{0,0}=\fmax(0)$}}
  \put( 9,-0.15){\line(0,1){0.3}}
  \put( 9,-0.7){\makebox(0,0){\footnotesize $g_{1,0}$}}
  \put(12,-0.15){\line(0,1){0.3}}
  \put(12,-0.7){\makebox(0,0){\footnotesize $f_{1,0}$}}
  \put(15,-0.15){\line(0,1){0.3}}
  \put(15,-0.7){\makebox(0,0){\footnotesize $g_{1,1}$}}
  \put(18,-0.3){\line(0,1){0.6}}
  \put(18,-1){\makebox(0,0){\footnotesize $f_{1,1}=\fmax(1)$}}
  \put(20,-0.5){\makebox(0,0){\footnotesize \dots}}
\end{picture}}
We require (for all $n,l\in\omega$)
\begin{itemize}
  \item $f_{n,l}\geq \Psi(n,g_{n,l},r(n))$ and
  \item $g_{n,l}\geq \varphi(\mathord<n)f_{n,l-1}^{n}$.
\end{itemize}
(Compare this with~\ref{asm:sufficientlydifferent}.)
Again it is clear that we can construct such sequences by induction.

Let $\myinvset$ be the set of $\nu: \omega\to\omega$ such that
$\nu(m)\leq m$ for all $m$. 
For  $\nu\in \myinvset$, we can define
$f_\nu:\omega \to \omega$ by $f_\nu(m)=f_{m,\nu(m)}$, and the same for $g_\nu$.
So we assign to each $\nu\in \myinvset$ cardinal characteristics
$\mycfa_{f_\nu,g_\nu}$ and $\myc_{f_\nu,g_\nu}$.

Assume that $X\subset \myinvset$ is countable such that
\begin{equation}\label{eq:wqr2}
  \text{for $\nu\neq \nu'$ in $X$ there is an $n(\nu, \nu')$
  such that $\nu(m)\neq \nu'(m)$ for all $m>n(\nu, \nu')$.}
\end{equation}
Then $(f_\nu,g_\nu)_{\nu\in X}$ is a suitable sequence 
as in Assumption~\ref{asm:sufficientlydifferent}.

\begin{uRem}
  We can of course find an uncountable set $X$ satisfying~\eqref{eq:wqr2}
  as well.
  We could try to define a forcing $P_I$ just as in the
  countable case, to force an uncountable version of
  Theorem~\ref{thm:ctbl}. However,
  we need ``separated support''~\ref{cor:variousbigness}(2)
  for~\eqref{eq:qwrr}. So we have to add appropriate
  requirements for conditions in $P$,
  in the style of $\trunklg^\text{min}$, this time depending 
  on the pair $\nu,\nu'$, to guarantee that 
  the maximum of the trunk-lengths at $\alpha\in I_\nu$
  and $\beta\in I_{\nu'}$ is bigger than the $n(\nu,\nu')$.
  However, such requirements lead to
  the following problem:
  Assume that $Y\subseteq \myinvset$ is countable and dense,
  and the domain of $p$ contains elements of $I_\nu$ 
  for each $\nu\in Y$. Then we cannot enlarge the domain 
  of $p$ to contain some $\nu'\notin Y$.\footnote{%
    In more detail: Let $f:Y\to\omega$ be such that 
    for all $\nu\in Y$, there is an $\alpha\in\dom(p)\cap I_\nu$
    such that $\trunklg(p,\alpha)+1<f(\nu)$. 
    Enumerate $Y$ as $\{\nu_0,\nu_1,\dots\}$.
    Then construct $\nu'\in\myinvset\setminus Y$
    the following way: Pick any $\nu^0\in Y$ and pick a finite
    $\nu'^0$ extending $\nu^0\restriction f(\nu^0)$, such that
    $\nu'^0(m)\neq \nu_0(m)$ for some $m$.
    Given $\nu'^l$, pick any $\nu^{l+1}\in Y$
    extending $\nu'^l$, and pick
    $\nu'^{l+1}$ extending $\nu^{l+1}\restriction f(\nu^{l+1})$
    such that $\nu'^{l+1}(m)\neq \nu_{l+1}(m)$ for some $m$.
    Set $\nu'=\bigcup_{l\in\omega} \nu'^l$. 
    Assume that there is a $q\leq p$ such that $\beta\in\dom(q)\cap I_{\nu'}$
    and $\trunklg(q,\beta)=m$. Only finitely many trunk-lengths
    in $\dom(p)$ were increased,
    so pick some $l$ such that $f(\nu^l)>m$ and 
    such that not trunk in $I_{\nu^l}$ was increased.
    By the definition of $f$,
    $\alpha\in\supp(q,m)$ for some $\alpha\in I_{\nu^{l}}$.
    $\nu'$ extends $\nu^{l}\restriction f(\nu^{l})$,
    so $\nu^{l}(m)=\nu'(m)$ (and $\nu^{l}\neq \nu'$), 
    which contradicts separated support.
  }
  So $p$ forces that the generic object does not contain
  anything in $I_{\nu'}$. But then our proofs do not
  work any more, cf.\ e.g.~\eqref{eq:bdomq}.
  To fix this problem, we will modify the forcing $P$ in the following
  way: As before, we choose for each $\epsilon\in\om1$ a cardinal
  $\kappa_\epsilon$ and the index set $I_\epsilon$ of
  size $\kappa_\epsilon$. However, we do not fix a $\nu\in\myinvset$
  for $\epsilon$. Instead, each condition $p$ chooses 
  $\nu(\epsilon)$ for each $\epsilon$ in its domain. 
  This makes Theorem~\ref{thm:uncountable} slightly weaker than
  Theorem~\ref{thm:ctbl}, since we do not know in the ground model
  which $\nu$ will be assigned to a $\kappa_\epsilon$.
\end{uRem}

We can now reformulate Theorem~\ref{thm:uncountable}:
\begin{Thm}
  Assume CH, assume
  that $\kappa_\epsilon=\kappa_\epsilon^\al0$ for $\epsilon\in\om1$,
  and that $(f_\nu,g_\nu)_{\nu\in\myinvset}$ are as above. Then there
  is a proper, $\al2$-cc, $\omega^\omega$-bounding partial order $R$ which
  forces:
  For each $\epsilon\in\om1$ there is a $\nu\in \myinvset$
  such that $\mycfa_{f_{\nu},g_{\nu}}=
  \myc_{f_{\nu},g_{\nu}}=\kappa_\epsilon$ 
  for all $\epsilon\in\om1$.
\end{Thm}
(Here $\myinvset$ denotes the set in $V$, not the evaluation of the definition
of $\myinvset$ in $V[G]$.)

As in the proof of Section~\ref{sec:countable},
we pick for each
$\nu\in \myinvset$ and $n\in\omega$ a creating pair $(\cK_\nu(n),\cS_\nu(n))$,
with $\bfH_\nu=f_\nu$ and $\bfF_\nu(n)=b$,
which is $(g_\nu(n),r(n))$-big, $r(n)$-halving and $(n,r(n))$-decisive.

We let $I$ be the disjoint union of $I_\epsilon$ ($\epsilon\in
\om1$), each $I_\epsilon$ has size
$\kappa_\epsilon$.

From here on, we assume CH. We now define the forcing notion $R$:
\begin{Def}
  A condition $p$ in $R$ consists of a countable 
  subset $\dom(p)$ of $I$, of objects $p(\alpha,n)$ for $\alpha\in \dom(p),
  n\in\omega$, and of functions $\trunklg(p):\dom(p)\to \omega$  
  and $\myinv(p):\dom(p)\to \myinvset$
  satisfying the following ($\alpha,\beta\in\dom(p)$): 
  \begin{itemize}
    \item $\myinv(p,\alpha)=\myinv(p,\beta)$ iff
      $\alpha,\beta$ are in the same $I_\epsilon$.
    \item
      If $n<\trunklg(p,\alpha)$, then
      $p(\alpha,n)\in \bfH_{\myinv(p,\alpha)}(n)$.
      \\
      $\bigcup_{n<\trunklg(p)}p(\alpha,n)$ is called trunk of $p$ at $\alpha$.
    \item
      If $n\geq \trunklg(p,\alpha)$, then $p(\alpha,n)\in \cK_{\myinv(p,\alpha)}(n)$
      and $\nor(p(\alpha,n))>0$.
    \item
      $\card{\supp(p,n)}<n$ for all $n>0$.
    \item
      Moreover, $\lim_{n\rightarrow\infty}(\card{\supp(p,n)}/n)=0$.
    \item
      $\lim_{n\rightarrow\infty}(\min(\{\nor(p(\alpha,n)):\, \alpha\in\supp(p,n)\}))= \infty$.
    \item \em{(Separated support.)}
      If $\alpha,\beta\in\supp(p,n)$, $\alpha\in I_\epsilon$,
      $\beta\in I_{\epsilon'}$,  and $\epsilon\neq \epsilon'$, then 
      $\myinv(p,\alpha)(n)\neq\myinv(p,\beta)(n)$. 
  \end{itemize}
\end{Def}
($\supp(p,n)$ is again the set of $\alpha\in I$ such that
$\trunklg(p,\alpha)\leq n$ .) 

Another way to formulate the last point is: If $\alpha,\beta\in\dom(p)$,
$\alpha\in I_\epsilon$, $\beta\in I_{\epsilon'}$, and $\epsilon\neq \epsilon'$,
then $\myinv(p,\alpha)$ and $\myinv(p,\beta)$ differ above some
$n(\myinv(p,\alpha),\myinv(p,\beta))$ as in~\eqref{eq:wqr2}, and
\begin{equation*}
  \max(\trunklg(p,\alpha),\trunklg(p,\beta))>n(\myinv(p,\alpha),\myinv(p,\beta)).
\end{equation*}

The order on $R$ is the natural modification of the one on $P$:
\begin{Def}
  For $p,q$ in $R$, $q\leq p$ if
  \begin{itemize}
    \item $\dom(q)\supseteq \dom(p)$,
    \item $\myinv(q,\alpha)=\myinv(p,\alpha)$ for all $\alpha\in\dom(p)$,
    \item if $\alpha\in\dom(p)$ and $n\in\omega$,
      then $q(\alpha,n)\in\cS_{\myinv(p,\alpha)}(p(\alpha,n))$,
    \item $\trunklg(q,\alpha)= \trunklg(p,\alpha)$ for
          all but finitely many $\alpha\in\dom(p)$.
  \end{itemize}
\end{Def}

$I_\epsilon$ is not a complete subforcing any more (conditions with disjoint
domains are generally not compatible, since the union can violate separated
support). But we still get:
\begin{Lem}
  $R$ is $\al2$-cc.
\end{Lem}
\begin{proof}
  Assume towards a contradiction that $A$ is an antichain
  of size $\al2$. By a $\Delta$-system argument, we 
  can assume that $\dom(p)\cap \dom(q)=u$
  for all distinct $p,q$ in $A$. We fix an enumeration 
  $\alpha^p_0,\alpha^p_1,\dots$ of $\dom(p)$ for each $p\in A$.
  By a pigeon hole argument, we 
  can assume that the following objects and statements
  are independent of $p\in A$ ($n,i\in\omega$, $\beta\in u$,
  $\epsilon\in\om1$): ``$\alpha^p_i=\beta$'', ``$\alpha^p_i\in I_\epsilon$'',
  $\trunklg(p,\alpha^p_i)$, $\myinv(p,\alpha^p_i)$, and $p(\alpha^p_i,n)$.

  So given distinct elements $p,q$ of $A$, we again 
  increase finitely many of the stems to guarantee that
  $\supp(p\cup q,n)$ has size less than $n$ for all $n$. 
  Then the resulting $r$ is a condition in $R$: To 
  see e.g.\ separated support, assume that 
  $\alpha,\beta\in\supp(r,n)$. We can assume that
  $\alpha=\alpha^p_i$ and $\beta=\alpha^q_j$ and that
  $\myinv(p,\alpha^p_i)\neq \myinv(q,\alpha^q_j)=\myinv(p,\alpha^p_j)$.
  Since $p$ satisfies separated support, 
  $\myinv(p,\alpha^p_i)(n)\neq \myinv(p,\alpha^p_j)(n)$.
\end{proof}

\begin{Lem} $R$ adds a generic real $\n\eta_\alpha$ for all $\alpha\in I$.
  In other words, the set of conditions $q$ with $\alpha\in\dom(q)$ 
  is dense.
\end{Lem}

\begin{proof} Assume $\alpha\in I_\epsilon$.
  Fix $p\in R$. We find a $q\leq p$ with $\dom(q)=\dom(p)\cup\{\alpha\}$.

  {\em Case 1:} $I_\epsilon\cap \dom(p)\neq \emptyset$.
  Then we pick $\beta\in I_\epsilon\cap \dom(p)$
  and choose $\trunklg(q,\alpha)>\trunklg(p,\beta)$
  big enough to guarantee $|\supp(q,n)|<n$ for all $n$.
  Then we choose any $q(\alpha,n)$ with sufficient norm (e.g.\ $n$).

  {\em Case 2:} Otherwise we again fix $\trunklg(q,\alpha)$
  big enough to guarantee $|\supp(q,n)|<n$ for all $n$,
  and we have to find
  some $\myinv(q,\alpha)$ satisfying
  separated support (for this $\trunklg(q,\alpha)$).
  Since $|\supp(p,n)|<n$ for all $n$, we can 
  find a $\nu'\in\myinvset$ such that $\nu'(n)$ is not in
  $\{\nu(n):\, \nu=\myinv(p,\beta),\beta\in\supp(p,n)\}$
  for any $n$.
  Set $\myinv(q,\alpha)=\nu'$.
  Then we again choose any $q(\alpha,n)$ with sufficiently 
  increasing norms.
  $q$ satisfies separated support:
  Assume $\beta\in I_{\epsilon'}\cap \supp(p,n)$.
  Then $\myinv(p,\beta)(n)\neq \nu'(n)=\myinv(q,\alpha)(n)$.
\end{proof}

It turns out that the proofs of Theorems~\ref{thm:prod}
and~\ref{thm:rapid} still work without any change:
\begin{Lem}
  $R$ is proper and $\omega\ho$-bounding.
  If $\delta\in\omega$, $\n \nu(n)$ is a $P$-name and
  $p\in P$ forces that
  $\n \nu(n)<\EXP(\fmax(n-1),n,n\cdot \delta)$ for all $n$,
  then there is a $q\leq p$ which
  $\mathord<n$-decides $\n \nu(n)$ for all $n$ .
\end{Lem}

\begin{proof}
  We define $\leq_n$ just as in the proof of Theorem~\ref{thm:prod}.
  Fusion still works: If $q$ is the limit of $p_n$,
  and each $p_n$ satisfies separated support, then so does $q$.
  The proof of pure decision does not require any changes.
  
  For rapid reading, note that each $\cK_\nu(n)$ is $(\fmax(n-1),r(n))$-big.
  Again, the same proof still works without changes.
\end{proof}

We can define the $R$-name $\myinv(\epsilon)$ for $\epsilon\in\om1$ to be
$\myinv(p,\alpha)$ for any $p$ in the generic filter and $\alpha\in\dom(p)\cap
I_\epsilon$.  Then we define the $R$-name $f_\epsilon$ to be
$f_{\myinv(\epsilon)}$, and the same for $g_\epsilon$.

We again get all items of Corollary~\ref{cor:variousbigness}, and can show:
\begin{Lem} 
  $R$ forces $\mycfa_{f_\epsilon,g_\epsilon}\leq \kappa_\epsilon$ and
  $\kappa_\epsilon\leq \myc_{f_\epsilon,g_\epsilon}$.
\end{Lem}

\begin{proof}
  The proofs of Lemmas~\ref{lem:forall} and~\ref{lem:exists} still work,
  if we assume that 
  $p_0$ determines $\myinv(\epsilon)$.
\end{proof}

\bibliographystyle{amsplain}
\bibliography{872}

\end{document}

%% file: pq_new.pstex_t
\begin{picture}(0,0)%
\includegraphics{pq_new.eps}%
\end{picture}%
\setlength{\unitlength}{4144sp}%
\begingroup\makeatletter\ifx\SetFigFont\undefined%
\gdef\SetFigFont#1#2#3#4#5{%
  \reset@font\fontsize{#1}{#2pt}%
  \fontfamily{#3}\fontseries{#4}\fontshape{#5}%
  \selectfont}%
\fi\endgroup%
\begin{picture}(10062,5274)(1876,-7348)
\put(5671,-7081){\makebox(0,0)[b]{\smash{{\SetFigFont{20}{24.0}{\familydefault}{\mddefault}{\updefault}{\color[rgb]{0,0,0}$q$}%
}}}}
\put(4051,-4741){\makebox(0,0)[b]{\smash{{\SetFigFont{20}{24.0}{\familydefault}{\mddefault}{\updefault}{\color[rgb]{0,0,0}$\val(p(2))\ni q(2)$}%
}}}}
\put(4051,-3661){\makebox(0,0)[b]{\smash{{\SetFigFont{20}{24.0}{\familydefault}{\mddefault}{\updefault}{\color[rgb]{0,0,0}$\cS(p(2))\ni q(3)$}%
}}}}
\put(4051,-7081){\makebox(0,0)[b]{\smash{{\SetFigFont{20}{24.0}{\familydefault}{\mddefault}{\updefault}{\color[rgb]{0,0,0}$>$}%
}}}}
\put(4051,-5641){\makebox(0,0)[b]{\smash{{\SetFigFont{20}{24.0}{\familydefault}{\mddefault}{\updefault}{\color[rgb]{0,0,0}Id}%
}}}}
\put(4051,-6451){\makebox(0,0)[b]{\smash{{\SetFigFont{20}{24.0}{\familydefault}{\mddefault}{\updefault}{\color[rgb]{0,0,0}Id}%
}}}}
\put(9496,-7081){\makebox(0,0)[b]{\smash{{\SetFigFont{20}{24.0}{\familydefault}{\mddefault}{\updefault}{\color[rgb]{0,0,0}$p$}%
}}}}
\put(11521,-7081){\makebox(0,0)[b]{\smash{{\SetFigFont{20}{24.0}{\familydefault}{\mddefault}{\updefault}{\color[rgb]{0,0,0}$q$}%
}}}}
\put(11476,-3166){\makebox(0,0)[b]{\smash{{\SetFigFont{20}{24.0}{\familydefault}{\mddefault}{\updefault}{\color[rgb]{0,0,0}$>M$}%
}}}}
\put(10486,-7081){\makebox(0,0)[b]{\smash{{\SetFigFont{20}{24.0}{\familydefault}{\mddefault}{\updefault}{\color[rgb]{0,0,0}$>_M$}%
}}}}
\put(10441,-5461){\makebox(0,0)[b]{\smash{{\SetFigFont{20}{24.0}{\familydefault}{\mddefault}{\updefault}{\color[rgb]{0,0,0}$=$}%
}}}}
\put(1891,-5281){\makebox(0,0)[rb]{\smash{{\SetFigFont{20}{24.0}{\familydefault}{\mddefault}{\updefault}{\color[rgb]{0,0,0}$\bfF(2)$}%
}}}}
\put(1891,-4201){\makebox(0,0)[rb]{\smash{{\SetFigFont{20}{24.0}{\familydefault}{\mddefault}{\updefault}{\color[rgb]{0,0,0}$\bfF(3)$}%
}}}}
\put(1891,-3031){\makebox(0,0)[rb]{\smash{{\SetFigFont{20}{24.0}{\familydefault}{\mddefault}{\updefault}{\color[rgb]{0,0,0}$\bfF(4)$}%
}}}}
\put(2521,-7081){\makebox(0,0)[b]{\smash{{\SetFigFont{20}{24.0}{\familydefault}{\mddefault}{\updefault}{\color[rgb]{0,0,0}$p$}%
}}}}
\put(1891,-6001){\makebox(0,0)[rb]{\smash{{\SetFigFont{20}{24.0}{\familydefault}{\mddefault}{\updefault}{\color[rgb]{0,0,0}$\bfF(1)$}%
}}}}
\put(1891,-6496){\makebox(0,0)[rb]{\smash{{\SetFigFont{20}{24.0}{\familydefault}{\mddefault}{\updefault}{\color[rgb]{0,0,0}$0$}%
}}}}
\put(1891,-3661){\makebox(0,0)[rb]{\smash{{\SetFigFont{20}{24.0}{\familydefault}{\mddefault}{\updefault}{\color[rgb]{0,0,0}$3$}%
}}}}
\put(1891,-4786){\makebox(0,0)[rb]{\smash{{\SetFigFont{20}{24.0}{\familydefault}{\mddefault}{\updefault}{\color[rgb]{0,0,0}$2$}%
}}}}
\put(1891,-5686){\makebox(0,0)[rb]{\smash{{\SetFigFont{20}{24.0}{\familydefault}{\mddefault}{\updefault}{\color[rgb]{0,0,0}$1$}%
}}}}
\put(1891,-7171){\makebox(0,0)[rb]{\smash{{\SetFigFont{29}{34.8}{\familydefault}{\mddefault}{\updefault}{\color[rgb]{0,0,0}$(a)$}%
}}}}
\put(8866,-7171){\makebox(0,0)[rb]{\smash{{\SetFigFont{29}{34.8}{\familydefault}{\mddefault}{\updefault}{\color[rgb]{0,0,0}$(b)$}%
}}}}
\put(8911,-3841){\makebox(0,0)[rb]{\smash{{\SetFigFont{20}{24.0}{\familydefault}{\mddefault}{\updefault}{\color[rgb]{0,0,0}$h$}%
}}}}
\end{picture}%

%% file: puredecbasic.pstex_t
\begin{picture}(0,0)%
\includegraphics{puredecbasic.eps}%
\end{picture}%
\setlength{\unitlength}{4144sp}%
\begingroup\makeatletter\ifx\SetFigFont\undefined%
\gdef\SetFigFont#1#2#3#4#5{%
  \reset@font\fontsize{#1}{#2pt}%
  \fontfamily{#3}\fontseries{#4}\fontshape{#5}%
  \selectfont}%
\fi\endgroup%
\begin{picture}(10467,5188)(1921,-7262)
\put(6076,-3661){\makebox(0,0)[rb]{\smash{{\SetFigFont{20}{24.0}{\familydefault}{\mddefault}{\updefault}{\color[rgb]{0,0,0}dec/}%
}}}}
\put(6076,-3976){\makebox(0,0)[rb]{\smash{{\SetFigFont{20}{24.0}{\familydefault}{\mddefault}{\updefault}{\color[rgb]{0,0,0}half}%
}}}}
\put(8776,-3661){\makebox(0,0)[rb]{\smash{{\SetFigFont{20}{24.0}{\familydefault}{\mddefault}{\updefault}{\color[rgb]{0,0,0}dec/}%
}}}}
\put(8776,-3976){\makebox(0,0)[rb]{\smash{{\SetFigFont{20}{24.0}{\familydefault}{\mddefault}{\updefault}{\color[rgb]{0,0,0}half}%
}}}}
\put(3826,-3661){\makebox(0,0)[rb]{\smash{{\SetFigFont{20}{24.0}{\familydefault}{\mddefault}{\updefault}{\color[rgb]{0,0,0}dec/}%
}}}}
\put(3826,-3976){\makebox(0,0)[rb]{\smash{{\SetFigFont{20}{24.0}{\familydefault}{\mddefault}{\updefault}{\color[rgb]{0,0,0}half}%
}}}}
\put(2476,-7126){\makebox(0,0)[b]{\smash{{\SetFigFont{20}{24.0}{\familydefault}{\mddefault}{\updefault}{\color[rgb]{0,0,0}$p$}%
}}}}
\put(4726,-7126){\makebox(0,0)[b]{\smash{{\SetFigFont{20}{24.0}{\familydefault}{\mddefault}{\updefault}{\color[rgb]{0,0,0}$p^0$}%
}}}}
\put(4726,-4741){\makebox(0,0)[b]{\smash{{\SetFigFont{20}{24.0}{\familydefault}{\mddefault}{\updefault}{\color[rgb]{0,0,0}\mydecr}%
}}}}
\put(8326,-7126){\makebox(0,0)[b]{\smash{{\SetFigFont{20}{24.0}{\familydefault}{\mddefault}{\updefault}{\color[rgb]{0,0,0}$\cdots$}%
}}}}
\put(9676,-7126){\makebox(0,0)[b]{\smash{{\SetFigFont{20}{24.0}{\familydefault}{\mddefault}{\updefault}{\color[rgb]{0,0,0}$p^{l-1}$}%
}}}}
\put(6976,-7126){\makebox(0,0)[b]{\smash{{\SetFigFont{20}{24.0}{\familydefault}{\mddefault}{\updefault}{\color[rgb]{0,0,0}$p^1$}%
}}}}
\put(6976,-2851){\makebox(0,0)[b]{\smash{{\SetFigFont{20}{24.0}{\familydefault}{\mddefault}{\updefault}{\color[rgb]{0,0,0}$\mathord\geq M$}%
}}}}
\put(9676,-2626){\makebox(0,0)[b]{\smash{{\SetFigFont{20}{24.0}{\familydefault}{\mddefault}{\updefault}{\color[rgb]{0,0,0}$\mathord\geq M$}%
}}}}
\put(6976,-4291){\makebox(0,0)[b]{\smash{{\SetFigFont{20}{24.0}{\familydefault}{\mddefault}{\updefault}{\color[rgb]{0,0,0}\mydecr}%
}}}}
\put(9676,-4111){\makebox(0,0)[b]{\smash{{\SetFigFont{20}{24.0}{\familydefault}{\mddefault}{\updefault}{\color[rgb]{0,0,0}\mydecr}%
}}}}
\put(4726,-3301){\makebox(0,0)[b]{\smash{{\SetFigFont{20}{24.0}{\familydefault}{\mddefault}{\updefault}{\color[rgb]{0,0,0}$\mathord\geq M$}%
}}}}
\put(11926,-7126){\makebox(0,0)[b]{\smash{{\SetFigFont{20}{24.0}{\familydefault}{\mddefault}{\updefault}{\color[rgb]{0,0,0}$S(p,M)$}%
}}}}
\put(11116,-5731){\makebox(0,0)[rb]{\smash{{\SetFigFont{20}{24.0}{\familydefault}{\mddefault}{\updefault}{\color[rgb]{0,0,0}big}%
}}}}
\put(1936,-5731){\makebox(0,0)[rb]{\smash{{\SetFigFont{20}{24.0}{\familydefault}{\mddefault}{\updefault}{\color[rgb]{0,0,0}$n$}%
}}}}
\end{picture}%

%% file: prod_pgtq_new.pstex_t
\begin{picture}(0,0)%
\includegraphics{prod_pgtq_new.eps}%
\end{picture}%
\setlength{\unitlength}{4144sp}%
\begingroup\makeatletter\ifx\SetFigFont\undefined%
\gdef\SetFigFont#1#2#3#4#5{%
  \reset@font\fontsize{#1}{#2pt}%
  \fontfamily{#3}\fontseries{#4}\fontshape{#5}%
  \selectfont}%
\fi\endgroup%
\begin{picture}(21492,7854)(-41729,-7078)
\put(-31049,389){\makebox(0,0)[b]{\smash{{\SetFigFont{29}{34.8}{\familydefault}{\mddefault}{\updefault}{\color[rgb]{0,0,0}$\omega$}%
}}}}
\put(-35549,389){\makebox(0,0)[b]{\smash{{\SetFigFont{29}{34.8}{\familydefault}{\mddefault}{\updefault}{\color[rgb]{0,0,0}$\omega$}%
}}}}
\put(-25469,389){\makebox(0,0)[b]{\smash{{\SetFigFont{29}{34.8}{\familydefault}{\mddefault}{\updefault}{\color[rgb]{0,0,0}$\omega$}%
}}}}
\put(-20294,-6451){\makebox(0,0)[b]{\smash{{\SetFigFont{29}{34.8}{\familydefault}{\mddefault}{\updefault}{\color[rgb]{0,0,0}$I$}%
}}}}
\put(-36539,-6451){\makebox(0,0)[b]{\smash{{\SetFigFont{29}{34.8}{\familydefault}{\mddefault}{\updefault}{\color[rgb]{0,0,0}$I$}%
}}}}
\put(-25919,-6451){\makebox(0,0)[b]{\smash{{\SetFigFont{29}{34.8}{\familydefault}{\mddefault}{\updefault}{\color[rgb]{0,0,0}$I$}%
}}}}
\put(-41669,389){\makebox(0,0)[b]{\smash{{\SetFigFont{29}{34.8}{\familydefault}{\mddefault}{\updefault}{\color[rgb]{0,0,0}$\omega$}%
}}}}
\put(-31724,-6451){\makebox(0,0)[b]{\smash{{\SetFigFont{29}{34.8}{\familydefault}{\mddefault}{\updefault}{\color[rgb]{0,0,0}$I$}%
}}}}
\put(-33659,-6901){\makebox(0,0)[b]{\smash{{\SetFigFont{34}{40.8}{\familydefault}{\mddefault}{\updefault}{\color[rgb]{0,0,0}$p$}%
}}}}
\put(-27989,-6901){\makebox(0,0)[b]{\smash{{\SetFigFont{34}{40.8}{\familydefault}{\mddefault}{\updefault}{\color[rgb]{0,0,0}$q$}%
}}}}
\put(-22364,-6901){\makebox(0,0)[b]{\smash{{\SetFigFont{34}{40.8}{\familydefault}{\mddefault}{\updefault}{\color[rgb]{0,0,0}$r$}%
}}}}
\put(-38609,-6901){\makebox(0,0)[b]{\smash{{\SetFigFont{34}{40.8}{\familydefault}{\mddefault}{\updefault}{\color[rgb]{0,0,0}$s$}%
}}}}
\put(-33974,-2086){\makebox(0,0)[b]{\smash{{\SetFigFont{29}{34.8}{\familydefault}{\mddefault}{\updefault}{\color[rgb]{0,0,0}$\supp$}%
}}}}
\put(-33614,-5371){\makebox(0,0)[b]{\smash{{\SetFigFont{29}{34.8}{\familydefault}{\mddefault}{\updefault}{\color[rgb]{0,0,0}$\text{trunk}$}%
}}}}
\put(-40319,-3256){\makebox(0,0)[b]{\smash{{\SetFigFont{29}{34.8}{\familydefault}{\mddefault}{\updefault}{\color[rgb]{0,0,0}$=p$}%
}}}}
\put(-28979,-5371){\makebox(0,0)[b]{\smash{{\SetFigFont{29}{34.8}{\familydefault}{\mddefault}{\updefault}{\color[rgb]{0,0,0}$\text{trunk}$}%
}}}}
\put(-26639,-5371){\makebox(0,0)[b]{\smash{{\SetFigFont{29}{34.8}{\familydefault}{\mddefault}{\updefault}{\color[rgb]{0,0,0}$\text{trunk}$}%
}}}}
\put(-26864,-2086){\makebox(0,0)[b]{\smash{{\SetFigFont{29}{34.8}{\familydefault}{\mddefault}{\updefault}{\color[rgb]{0,0,0}$\supp$}%
}}}}
\put(-29474,-1366){\makebox(0,0)[b]{\smash{{\SetFigFont{29}{34.8}{\familydefault}{\mddefault}{\updefault}{\color[rgb]{0,0,0}$\supp$}%
}}}}
\put(-37529,-961){\makebox(0,0)[b]{\smash{{\SetFigFont{20}{24.0}{\familydefault}{\mddefault}{\updefault}{\color[rgb]{0,0,0}$\nor\mathord>M$}%
}}}}
\put(-40094,-961){\makebox(0,0)[b]{\smash{{\SetFigFont{20}{24.0}{\familydefault}{\mddefault}{\updefault}{\color[rgb]{0,0,0}$\nor\mathord>M$}%
}}}}
\put(-41714,-2491){\makebox(0,0)[rb]{\smash{{\SetFigFont{29}{34.8}{\familydefault}{\mddefault}{\updefault}{\color[rgb]{0,0,0}$h$}%
}}}}
\put(-39419,-2986){\makebox(0,0)[lb]{\smash{{\SetFigFont{34}{40.8}{\familydefault}{\mddefault}{\updefault}{\color[rgb]{0,0,0}$=$}%
}}}}
\put(-28529,-1366){\makebox(0,0)[lb]{\smash{{\SetFigFont{34}{40.8}{\familydefault}{\mddefault}{\updefault}{\color[rgb]{0,0,0}$=$}%
}}}}
\put(-21194,-466){\makebox(0,0)[b]{\smash{{\SetFigFont{20}{24.0}{\familydefault}{\mddefault}{\updefault}{\color[rgb]{0,0,0}$\nor\mathord>M$}%
}}}}
\put(-24119,-2671){\makebox(0,0)[b]{\smash{{\SetFigFont{29}{34.8}{\familydefault}{\mddefault}{\updefault}{\color[rgb]{0,0,0}$=q$}%
}}}}
\put(-21194,-1051){\makebox(0,0)[b]{\smash{{\SetFigFont{29}{34.8}{\familydefault}{\mddefault}{\updefault}{\color[rgb]{0,0,0}$=q$}%
}}}}
\put(-24794,-1861){\makebox(0,0)[rb]{\smash{{\SetFigFont{29}{34.8}{\familydefault}{\mddefault}{\updefault}{\color[rgb]{0,0,0}$h$}%
}}}}
\put(-22724,-3886){\makebox(0,0)[b]{\smash{{\SetFigFont{20}{24.0}{\familydefault}{\mddefault}{\updefault}{\color[rgb]{0,0,0}$|\supp|/h<1/(M+1)$}%
}}}}
\end{picture}%

%% file: motivation.pstex_t
\begin{picture}(0,0)%
\includegraphics{motivation.eps}%
\end{picture}%
\setlength{\unitlength}{4144sp}%
\begingroup\makeatletter\ifx\SetFigFont\undefined%
\gdef\SetFigFont#1#2#3#4#5{%
  \reset@font\fontsize{#1}{#2pt}%
  \fontfamily{#3}\fontseries{#4}\fontshape{#5}%
  \selectfont}%
\fi\endgroup%
\begin{picture}(20638,8026)(-42675,-7226)
\put(-35549,389){\makebox(0,0)[b]{\smash{{\SetFigFont{29}{34.8}{\familydefault}{\mddefault}{\updefault}{\color[rgb]{0,0,0}$\omega$}%
}}}}
\put(-35639,-2491){\makebox(0,0)[rb]{\smash{{\SetFigFont{34}{40.8}{\familydefault}{\mddefault}{\updefault}{\color[rgb]{0,0,0}$n$}%
}}}}
\put(-29429,-6451){\makebox(0,0)[b]{\smash{{\SetFigFont{29}{34.8}{\familydefault}{\mddefault}{\updefault}{\color[rgb]{0,0,0}$I$}%
}}}}
\put(-32624,-6676){\makebox(0,0)[b]{\smash{{\SetFigFont{29}{34.8}{\familydefault}{\mddefault}{\updefault}{\color[rgb]{0,0,0}$S(p,M)$}%
}}}}
\put(-24929,-6676){\makebox(0,0)[b]{\smash{{\SetFigFont{29}{34.8}{\familydefault}{\mddefault}{\updefault}{\color[rgb]{0,0,0}$q'\leq q^\omega$}%
}}}}
\put(-27899,389){\makebox(0,0)[b]{\smash{{\SetFigFont{29}{34.8}{\familydefault}{\mddefault}{\updefault}{\color[rgb]{0,0,0}$\omega$}%
}}}}
\put(-42299,-6991){\makebox(0,0)[b]{\smash{{\SetFigFont{41}{49.2}{\familydefault}{\mddefault}{\updefault}{\color[rgb]{0,0,0}(a)}%
}}}}
\put(-35549,-6991){\makebox(0,0)[b]{\smash{{\SetFigFont{41}{49.2}{\familydefault}{\mddefault}{\updefault}{\color[rgb]{0,0,0}(b)}%
}}}}
\put(-27854,-6991){\makebox(0,0)[b]{\smash{{\SetFigFont{41}{49.2}{\familydefault}{\mddefault}{\updefault}{\color[rgb]{0,0,0}(c)}%
}}}}
\put(-42299,389){\makebox(0,0)[b]{\smash{{\SetFigFont{29}{34.8}{\familydefault}{\mddefault}{\updefault}{\color[rgb]{0,0,0}$\omega$}%
}}}}
\put(-37619,-6451){\makebox(0,0)[b]{\smash{{\SetFigFont{29}{34.8}{\familydefault}{\mddefault}{\updefault}{\color[rgb]{0,0,0}$I$}%
}}}}
\put(-40049,-6991){\makebox(0,0)[b]{\smash{{\SetFigFont{34}{40.8}{\familydefault}{\mddefault}{\updefault}{\color[rgb]{0,0,0}$\dom(p)$}%
}}}}
\put(-42389,-3841){\makebox(0,0)[rb]{\smash{{\SetFigFont{34}{40.8}{\familydefault}{\mddefault}{\updefault}{\color[rgb]{0,0,0}$n$}%
}}}}
\put(-23039,-1816){\makebox(0,0)[b]{\smash{{\SetFigFont{29}{34.8}{\familydefault}{\mddefault}{\updefault}{\color[rgb]{0,0,0}$q^1$}%
}}}}
\put(-26684,-4246){\makebox(0,0)[b]{\smash{{\SetFigFont{29}{34.8}{\familydefault}{\mddefault}{\updefault}{\color[rgb]{0,0,0}$q^1$}%
}}}}
\put(-26414,-3391){\makebox(0,0)[b]{\smash{{\SetFigFont{29}{34.8}{\familydefault}{\mddefault}{\updefault}{\color[rgb]{0,0,0}$q^2$}%
}}}}
\put(-23579,-1186){\makebox(0,0)[b]{\smash{{\SetFigFont{29}{34.8}{\familydefault}{\mddefault}{\updefault}{\color[rgb]{0,0,0}$q^2$}%
}}}}
\put(-26144,-2491){\makebox(0,0)[b]{\smash{{\SetFigFont{29}{34.8}{\familydefault}{\mddefault}{\updefault}{\color[rgb]{0,0,0}$q^3$}%
}}}}
\put(-23984,-601){\makebox(0,0)[b]{\smash{{\SetFigFont{29}{34.8}{\familydefault}{\mddefault}{\updefault}{\color[rgb]{0,0,0}$q^3$}%
}}}}
\put(-26999,-736){\makebox(0,0)[b]{\smash{{\SetFigFont{29}{34.8}{\familydefault}{\mddefault}{\updefault}{\color[rgb]{0,0,0}$q'$}%
}}}}
\put(-24749,-736){\makebox(0,0)[b]{\smash{{\SetFigFont{29}{34.8}{\familydefault}{\mddefault}{\updefault}{\color[rgb]{0,0,0}$q^\omega$}%
}}}}
\put(-40049,-5371){\makebox(0,0)[b]{\smash{{\SetFigFont{29}{34.8}{\familydefault}{\mddefault}{\updefault}{\color[rgb]{0,0,0}$\text{trunk}$}%
}}}}
\put(-34199,-3796){\makebox(0,0)[b]{\smash{{\SetFigFont{29}{34.8}{\familydefault}{\mddefault}{\updefault}{\color[rgb]{0,0,0}$\text{big}$}%
}}}}
\put(-40949,-4471){\makebox(0,0)[b]{\smash{{\SetFigFont{29}{34.8}{\familydefault}{\mddefault}{\updefault}{\color[rgb]{0,0,0}$\supp(p,n)$}%
}}}}
\put(-34199,-4471){\makebox(0,0)[b]{\smash{{\SetFigFont{29}{34.8}{\familydefault}{\mddefault}{\updefault}{\color[rgb]{0,0,0}$\supp(p,n)$}%
}}}}
\end{picture}%